\RequirePackage{fix-cm}
\documentclass[smallextended]{svjour3}       
\smartqed  

\usepackage{amsmath}  
\usepackage{fancyhdr}
\usepackage{extramarks}
\usepackage{amsfonts}
\usepackage{tikz}
\usepackage{algorithm}
\usepackage[noend]{algpseudocode}
\usepackage{amssymb}
\usepackage{mathrsfs}
\usepackage{algorithmicx}
\usepackage{stmaryrd}
\usepackage{color}
\usepackage{verbatim}
\usepackage{hyperref}
\hypersetup{hidelinks}
\usepackage{caption}
\usepackage{graphicx}
\usepackage{float} 
\usepackage{booktabs}
\usepackage{subcaption}
\usepackage{cases,mathdots,colordvi,appendix,url,multirow,lscape}
\usepackage{geometry}
\usepackage{caption}
\usepackage[justification=centering]{caption}

\makeatletter
\@addtoreset{equation}{section}
\makeatother

\definecolor{BLUE}{rgb}{0,0,1}

\DeclareMathOperator*{\argmin}{argmin}

\begin{document}

\title{On the robust isolated calmness of a class of nonsmooth optimizations on Riemannian manifolds and its applications}
\author{Chenglong Bao \and
	Chao Ding   \and
	Yuexin Zhou        
}
\institute{
           C.L. Bao \at
              Yau Mathematical Sciences Center, Tsinghua University, P.R. China and Yanqi Lake Beijing Institute of Mathematical Sciences and Applications, P.R. China. C. Bao is supported by the National Key R\&D Program of China (No.~2021YFA1001300) and the National Natural Science Foundation of China (No.~12271291).\\ 
              \email{clbao@mail.tsinghua.edu.cn}
              \and
              C. Ding \at
              State Key Laboratory of Mathematical Sciences, Academy of Mathematics and Systems Science, Chinese Academy of Sciences, Beijing 100190, P.R. China; School of Mathematical Sciences, University of Chinese Academy of Sciences, Beijing 100049, P.R. China; Institute of Applied Mathematics, Academy of Mathematics and Systems Science, Chinese Academy of Sciences, Beijing 100190, P.R. China. The research of this author was supported in part by the National Key R\&D Program of China (2021YFA1000300, 2021YFA1000301), National Natural Science Foundation of China (12071464) and CAS Project for Young Scientists in Basic Research (YSBR-034). \\ 
              \email{dingchao@amss.ac.cn}
              \and
              Y.X. Zhou \at
              Institute of Applied Mathematics, Academy of Mathematics and Systems Science, Chinese Academy of Sciences, Beijing, P.R. China, School of Mathematical Sciences, University of Chinese Academy of Science, Beijing, P.R. China. \\
              \email{zhouyuexin19@mails.ucas.ac.cn}           
}

\date{This version: April 8, 2025}

\maketitle

\begin{abstract}
	This paper studies the robust isolated calmness property of the KKT solution mapping of a class of nonsmooth optimization problems on {\color{black}Riemannian manifolds}. The manifold versions of the Robinson constraint qualification, the strict Robinson constraint qualification, and the second order conditions are defined and discussed. We show that the robust isolated calmness of the  {\color{black} KKT solution mapping is equivalent to satisfying the M-SRCQ and  M-SOSC conditions.} Furthermore, under the above two conditions, we show that the Riemannian augmented Lagrangian method achieves a local linear convergence rate. Finally, we verify the proposed conditions and demonstrate the convergence rate on two minimization problems over the sphere and the manifold of fixed rank matrices.

\end{abstract}

\section{Introduction}\label{sec:intro}

In recent years, manifold optimization has become an important class of constrained optimization problems and has found various applications in many tasks, {\color{black}such as computer vision, signal processing, statistical learning, numerical linear algebra, and machine learning (cf. e.g., \cite{L05,R84,V13,ZHT06}). }For a comprehensive study of manifold optimization, see~\cite{AMS09,B20,HLWY20}. In this paper, we consider the following nonsmooth manifold optimization problem:
\begin{equation}\label{eq:prime}
	\begin{array}{ll}
		\min & f(x)+\theta(g_1(x)) \\
		\text { s.t. } 
		& g_2(x) \in \mathcal{Q},\\
		& x\in \mathcal{M},
	\end{array}
\end{equation}
where $\mathcal{M}$ is a {\color{black}finite-dimensional} smooth Riemannian manifold, $f:\mathcal{M} \rightarrow \mathbb{R}$, $g_1:\mathcal{M} \rightarrow \mathbb{Y}$ and $g_2:\mathcal{M} \rightarrow \mathbb{Z}$ are twice continuously differentiable functions {\color{black}from manifold $\mathcal{M}$ to vector spaces}, $\mathbb{Y}$ and $\mathbb{Z}$  are two Euclidean spaces each equipped with a scalar
product $\langle\cdot, \cdot\rangle$ and its induced norm $\|\cdot\|$,  $\theta : \mathbb{Y}\rightarrow \mathbb{R}$ is a proper closed convex function, and $\mathcal{Q} \subset \mathbb{Z}$ is a nonempty closed convex set. Many applications arising from emerging fields can be cast into the form  (\ref{eq:prime}), e.g., compressed modes \cite{OLCO13}, sparse principal component analysis \cite{ZHT06}, constrained sparse principal analysis \cite{LZ12} and robust matrix completion~\cite{CA16}, here we list two of these examples in the following.  See \cite{AH19} for more details.

{\color{black}
	\begin{example}{\textbf{Sparse principle component analysis(SPCA)}\cite{ZHT06}}\label{exam:spca}
	
	Let $A \in \mathbb{R}^{m \times n}$ be the given data matrix where $n$ is the number of variables and $m$ is the number of observations. Set $f(X)=-\operatorname{tr}(X^{\top}A^{\top}AX)$, $g(X)=X$,$\theta(\cdot)=\mu\|\cdot\|_{1}$, $\mathcal{M}=\operatorname{St}(n,r)$. Then the SPCA problem can be formulated as 
	\begin{equation}
		\begin{array}{ll}
			\min _{X \in \mathbb{R}^{n \times r}} & -\operatorname{tr}\left(X^{\top} A^{\top} A X\right)+\mu\|X\|_{1} \\
			\text { s.t. } & X^{\top} X=I_{r}.
		\end{array}
	\end{equation}
This model is a nonsmooth optimization problem on the Stiefel manifold, which aims to find the principal components while imposing the L1-norm to induce sparsity. It is crucial for identifying the most important features in high-dimensional data.
\end{example}

\begin{example}{\textbf{Robust matirx completion(RMC)}\cite{CA16}}\label{exam:rmc}
	We are now considering a basic robust matrix completion problem. 
	For a given $A\in \mathbb{R}^{m\times n}$, let $g(X)=P_{\Omega}(X-A)$ and $\theta(\cdot)=\mu\|\cdot\|_{1}$. Here,  $P_{\Omega}$ is the projector defined by $\left(P_{\Omega}(X)\right)_{i j}=X_{i j}$ if $(i, j) \in \Omega$ and $0$ otherwise. By setting $\mathcal{M}=Fr(m,n,r):=\{X\in \mathbb{R}^{m\times n}:rank(X)=r\}$, we obtain the following robust matrix completion problem
	\begin{equation}
		\begin{array}{ll}
			\min _{X \in \mathbb{R}^{m \times n}} & \left\|P_{\Omega}(X-A)\right\|_{1} \\
			\text { s.t. } & X \in Fr(m,n,r).
		\end{array}
	\end{equation}
	This model is a nonsmooth optimization problem on the fixed-rank manifold, which aims to recover a low-rank approximation of a matrix while minimizing the impact of outliers and noise. The robust matrix completion is essential for accurately reconstructing missing data in the presence of corruption or errors.
\end{example}
}

Many algorithms for solving nonsmooth optimization problems have been extended from the Euclidean space to Riemannian manifolds, such as the subgradient methods~\cite{FO98,GH16}, the proximal gradient methods~\cite{CMMZ20,HW19,HW21}, the alternating direction methods of multipliers (ADMM)~\cite{KGB16,LO14}, the proximal point methods~\cite{CDMS21,FO02}, and the augmented Lagrangian methods (ALM)~\cite{DP19,ZBDZ21}. However, the theoretical results for nonsmooth manifold optimization problems appear to be significantly fewer compared to those in Euclidean settings. To the best of our knowledge, perturbation analysis for problem~\eqref{eq:prime}, which is closely related to the convergence analysis of numerical algorithms, has not yet been established.

\begin{figure}[h]
	\centering
	\includegraphics[scale=0.5]{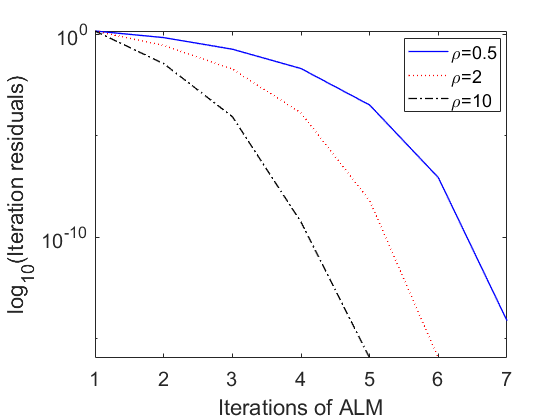}
	\caption{The Riemannian ALM for solving  (\ref{exam:sphere}) with different penalty parameters $\rho$}
	\label{fig:ex-1}
\end{figure}

One of the essential perturbation properties for optimization problems is the robust isolated calmness of the Karush-Kuhn-Tucker (KKT) solution mapping under perturbations \cite[Definition 2]{DSZ17}. The isolated calmness is a Lipschitz-like property of set-value mappings, which is first proposed by Robinson  \cite{R79} and generalized in \cite{D09} as follows: a mapping $F: \mathbb{R}^{m} \rightarrow \mathbb{R}^{n}$ is said to have the isolated calmness property if there exist a constant $\kappa \geq 0$ and neighborhoods $\mathcal{V}$ of $\bar{y}$ and $\mathcal{U}$ of $\bar{x}$ such that $
\|y-\bar{y}\| \leq \kappa\|x-\bar{x}\|$ when $y \in F(x) \cap \mathcal{V}$ and $x \in \mathcal{U}$. The robust version is defined in \cite{DSZ17} with an additional requirement: $F(x) \cap \mathcal{V}\neq \emptyset$ for all $x \in \mathcal{U}$. This property is crucial for establishing the linear convergence rate of numerical algorithms. For example, it is used to analyze the convergence of sequential quadratic programming (SQP) in \cite{B94} for nonlinear programs (NLPs), where the property is obtained under the strict Mangasarian-Fromovitz constraint qualification (SMFCQ) and the second order sufficient condition (SOSC). The recent work~\cite{DSZ17} successfully characterizes the robust isolated calmness for optimization problems with $\mathcal{C}^2$-cone reducible constraints  (cf. Definition \ref{def:c2reducible}) in Euclidean space. Specifically, under the Robinson constraint qualification (RCQ), it is proved that the robust isolated calmness of the KKT mapping holds at a local solution if and only if the strict Robinson constraint qualification (SRCQ) and the SOSC are satisfied at the given point. The main challenge for establishing the perturbation analysis for manifold optimization is to construct a proper perturbed problem since the canonical perturbation $\left\langle a,\cdot\right\rangle$ is no longer a linear function on manifolds. In this work, we first construct the perturbed manifold problem through a locally equivalent Euclidean problem using the normal coordinate chart around a KKT point and define the SRCQ and SOSC conditions on manifolds. Furthermore, we show that the robust isolated calmness of the KKT solution mapping of \eqref{eq:prime} is equivalent to the manifold SRCQ (M-SRCQ) condition and the manifold SOSC (M-SOSC) condition (the definitions are given in Section \ref{sec:ric}) hold at solution points, which is the manifold extension of \cite{DSZ17}.\\
\indent An important application of the robust isolated calmness of the KKT solution mapping is the local convergence analysis of the augmented Lagrangian method (ALM). The classical ALM is proposed by Hestenes \cite{H69} and Powell \cite{P69} for equality constraints and is extended to nonlinear programming by Rockafellar  \cite{R73}. The convergence analysis of ALM under Euclidean settings has been extensively studied for decades. The classical result of the local linear convergence rate of ALM for NLPs often requires the linear independence constraint qualification (LICQ) and the second order sufficient condition, e.g.\ \cite{B82,CGT00,WN99}. For conic programs such as the nonlinear second order cone programs and semidefinite programs, Liu and Zhang~\cite{LZ08}, and Sun et al.~\cite{SSZ08}, respectively, obtained the local linear convergence rate of ALM under the constraint non-degeneracy condition and the strong second order sufficient condition at KKT points{\color{black}. }Recently, for general nonlinear optimizations involving  $\mathcal{C}^2$-cone reducible constraints, Kanzow and Steck \cite{KS19} obtain the primal-dual linear convergence result under the assumption of the robust isolated calmness of the KKT solution mapping, which is equivalent to the SRCQ and SOSC conditions by \cite{DSZ17}. \\
\indent The Riemannian Augmented Lagrangian Method (ALM) for nonsmooth constrained manifold optimization was recently introduced in \cite{ZBDZ21} (see also \cite{DP19}, \cite{LB20}). {\color{black}In \cite{LB20}, the global convergence of ALM to a KKT point is established under the LICQ. \cite{ZBDZ21} also proves the convergence of the iteration sequence of the Riemannian ALM to a KKT point, assuming suitable constraint qualifications, and presents numerical experiments demonstrating its superior performance compared to existing methods. Additionally, \cite{YS22l} extends the approximate KKT (AKKT) conditions to nonlinear optimization on Riemannian manifolds and proposes an ALM that globally converges to points satisfying AKKT. \cite{ACFH24} formulates several well-known constraint qualifications from Euclidean optimization, ensuring the global convergence of the Riemannian ALM (RALM) without requiring the boundedness of the multiplier set.} Motivated by the promising results in~\cite{ZBDZ21}, one natural question is whether the local convergence rate of the Riemannian ALM for the nonsmooth manifold optimization \eqref{eq:prime}  can similarly be obtained as in \cite{KS19} under the robust isolated calmness for the KKT solution, or equivalently, the assumption of the M-SRCQ and M-SOSC conditions. First, we consider a simple example:
\begin{equation}\label{exam:sphere}  
	\begin{array}{ll}
		\min & x_2^2+|x_1-x_2| \\
		\text { s.t. } 
		& 2x_1+x_2\geq 0,\\
		& {\color{black}x\in \mathcal{M}=\mathcal{S}^1.}
	\end{array} 
\end{equation}
{\color{black}Here, $\mathcal{S}^1=\{x\in \mathbb{R}^2|x_1^2+x_2^2=1\}$ is the unit sphere in $ \mathbb{R}^2$.} {\color{black}It is verified that the unique KKT solution of \eqref{exam:sphere} is  $x^*=({\sqrt{2}}/{2},{\sqrt{2}}/{2})^T$. By computing the KKT conditions, we can obtain the corresponding multipliers $y^*={\sqrt{2}}/{2}$ for the nonsmooth term and $z^*=0$ for the inequality constraint.} By applying the Riemannian ALM (see Algorithm \ref{alg:alm} for details) and noting that the corresponding ALM subproblems can be solved exactly, it is observed from Figure \ref{fig:ex-1} that  the distance between the iteration $(x^k,y^k,z^k)$ and the solution $(x^*,y^*,z^*)$ converges linearly as
$k$ becomes sufficiently large.  Moreover, Figure \ref{fig:ex-1} indicates that the linear convergence rate decreases as the penalty parameter $\rho^{k}$ increases.
On the other hand, it will be verified in Section \ref{sec:ric} (Remark \ref{rm:example-sosc-srcq}) that for problem (\ref{exam:sphere}), the M-SRCQ holds at $x^*$ with respect to $(y^*,z^*)$ and the M-SOSC holds at $x^*$. Inspired by \cite{KS19}, in this paper, we show that the Riemannian ALM indeed has a local linear convergence rate under the M-SRCQ and M-SOSC assumptions. To the best of our knowledge, this is the first time the local convergence rate of ALM on manifolds has been obtained. We use two examples of nonsmooth optimization over a sphere and a fixed rank manifold to verify the M-SRCQ and M-SOSC conditions for different problems, and to illustrate the theoretical results obtained in Section \ref{sec:numerical}. We will also show that the M-SRCQ condition is satisfied at all KKT points for a class of nonsmooth optimization problems on the sphere.

When the manifold $\mathcal{M}$ in problem (\ref{eq:prime}) is embedded in a Euclidean space $\mathbb{X}$ (cf. \cite[Definition 3.10]{B20}),  it can be locally expressed as the equality constraint {\color{black}$\left\lbrace x\in \mathbb{X}\mid h(x)=0\right\rbrace$, where $h:\mathbb{X}\to \mathbb{R}^m$ for dimension $m$.} In applications, many manifold constraints for optimization are embedded, for example, the Euclidean space itself, the Stiefel manifold $\{X\in \mathbb{R}^{n\times p}\mid X^{\top}X=I_p\}$ and the fixed rank manifold $\{X\in\mathbb{R }^{m\times n}\mid \operatorname{rank}(X)=r\}$. In this case, the perturbation properties of problem (\ref{eq:prime}) can be studied through the classical approaches. {\color{black}If $h'(x)$ is onto $\mathbb{R}^m$, }we will show in Section \ref{sec:ric} (Remarks \ref{rm:em-srcq-equi} and \ref{rm:em-sosc-equi}) that the Euclidean SRCQ and the Euclidean SOSC conditions for problem (\ref{eq:prime})  are equivalent to the M-SRCQ and M-SOSC conditions. However, if $h'(x)$ is not of full row rank, we find that the Euclidean SRCQ condition cannot be satisfied even when the problem has no other constraints. Additionally, when $\mathcal{M}$ is no longer an embedded submanifold, applying traditional ALM to problem (\ref{eq:prime}) seems challenging. Based on these two observations, in this paper, we study the perturbation properties of problem~(\ref{eq:prime}) in the manifold setting and the local convergence rate of the Riemannian augmented Lagrangian method directly.

The rest of the paper is organized as follows. In Section \ref{sec:preli} we review some background of smooth manifolds and set-valued mappings. In Section \ref{sec:ric}, we define the constraint qualifications and second order optimality conditions on manifold, and characterize the robust isolated calmness of the KKT solution mapping for the perturbed problem introduced in this section. The Riemannian augmented Lagrangian method and its local convergence analysis are presented in Section \ref{sec:alm}. The applications and numerical results are contained in Section \ref{sec:numerical}. Finally, we make some remarks in Section \ref{sec:conslusion}.

\section{Preliminaries and notations }\label{sec:preli}

We first introduce some basic concepts of manifolds and set-valued mappings that will be used in our discussion. {\color{black}The notations that will be used in the following discussions are listed in Table \ref{table:notations}. }The properties of manifolds discussed here can be found in \cite{AMS09,K11,L13}. Information on set-valued mappings and directional derivatives can be found in \cite{BS13,D09}.

\renewcommand{\arraystretch}{1.3}
\begin{table}[h]
	\centering
	\caption{Notations used in this work}
	\label{table:notations}
	\begin{tabular}{|c | l | }
		\hline
	\textbf{	Notations} & \multicolumn{1}{c|}{\textbf{Descriptions}} \\
		\hline
	$\mathbb{X}$, $\mathbb{Y}$, $\mathbb{Z}$ &  Euclidean spaces\\
	\hline
		$\mathcal{M}$ & a $n$-dimensional smooth Riemannian manifold\\
	\hline
		$\varphi$, $(\mathcal{U},\varphi)$ & the chart of $\mathcal{M}$\\
		\hline
		$ \mathfrak{F}_{x}(\mathcal{M}) $ & the set of smooth real-valued functions defined around $ x \in \mathcal{M}$\\
		\hline
		 $ \xi_{x} $& a tangent vector at $x$\\
		 \hline
		 $ T_x\mathcal{M} $&  the tangent space at $x$\\
		 	 \hline
		 $ N_x\mathcal{M} $&  the normal space at $x$\\
		 	 \hline
		 $ T\mathcal{M} $&  the tangent bundle of $\mathcal{M}$\\
		\hline
	$\hat{x}$, $\hat{\xi}$, $\hat{f}$ & the coordinate representation of $x$, $\xi$, $f$\\
	\hline
$ \langle\cdot, \cdot\rangle_{x} $, 	$ \langle\cdot, \cdot\rangle $ & the Riemannian metric at $x\in \mathcal{M}$\\
	\hline
$ DF(x) $& the differential at $ x \in\mathcal{M}$ for Riemannian mapping $F$\\
\hline
$ h'(x) $& the differential at $ x \in\mathbb{X}$ for Euclidean mapping $h$\\
\hline
$\operatorname{grad}f(x)$ &the Riemannian gradient of $f$ at $x\in \mathcal{M}$\\
\hline
$\nabla h(x)$ & the Euclidean gradient of $h$ at $x\in \mathbb{X}$\\
\hline
$\operatorname{Hess}f(x)$ &the Riemannian Hessian of $f$ at $x\in \mathcal{M}$\\
\hline
$ \Pi_{x} $ & the projection onto $ T_{x} \mathcal{M} $\\
\hline
$	f^{\circ}(x ; v)$&the generalized directional derivative of $f$ at $x \in \mathcal{M}$ in  direction $v \in T_{x} \mathcal{M}$\\
\hline
$\partial f(x)$ & the Clarke subdifferential of $f$ at $x \in \mathcal{M}$\\
\hline
$ \operatorname{Exp}_{x}$ & the exponential map at $ x \in \mathcal{M} $\\
\hline
$\mathcal{T}_{\mathcal{K}}(y)$ & the tangent cone to convex set $\mathcal{K}$ at $y\in \mathbb{Y}$\\
\hline
$\mathcal{N}_{\mathcal{K}}(y)$ & the normal cone to convex set $\mathcal{K}$ at $y\in \mathbb{Y}$\\
\hline
$d(x,x^{\prime})$ & the Riemannian distance of $x,x^{\prime}\in\mathcal{M}$\\
\hline
$\operatorname{dist}(y,\mathcal{U})$ & the distance from  point $y$ to set $\mathcal{U}$\\
\hline
\end{tabular}
\end{table}

Given a smooth $ n $-dimensional manifold $ \mathcal{M} $ and a point $ x \in \mathcal{M}$, {\color{black}a bijection $\varphi$ of a subset $\mathcal{U}$ of $\mathcal{M}$ onto an open subset of $\mathbb{R}^n$ is called a $n$-dimensional chart of $\mathcal{M}$, denoted by $(\mathcal{U},\varphi)$. When no confusion arises, we simply use $\varphi$ instead of $(\mathcal{U},\varphi)$. Let $f$ be a function from manifold $\mathcal{M}$ into $\mathbb{R}^m$. Choosing $x\in\mathcal{M}$ and a chart $\varphi$ around $x$, function $\hat{f}= f\circ \varphi^{-1}$ is called a coordinate representation of $f$. We say that $f$ is differentiable or smooth at $x$ if $\hat{f}$ is differentiable or smooth at $\varphi(x)$.} Let $ \mathfrak{F}_{x}(\mathcal{M}) $ represent the set of all smooth real-valued functions defined on a neighborhood of $ x $. A mapping $ \xi_{x} $ from $ \mathfrak{F}_{x}(\mathcal{M}) $ to $ \mathbb{R} $ is defined such that there exists a curve $ \gamma $ on $ \mathcal{M} $ with $ \gamma(0) = x $ satisfying $ \xi_{x} f := \dot{\gamma}(0)f := \left.\frac{d(f(\gamma(t)))}{dt}\right|_{t=0} $ for any $ f \in \mathfrak{F}_{x}(\mathcal{M}) $. This mapping $ \xi_{x} $ is referred to as a tangent vector, and the set of all such tangent vectors at $ x $ forms the tangent space $ T_x\mathcal{M} $. When $ \mathcal{M} $ is embedded in a Euclidean space $ \mathbb{X} $, the normal space $ N_x\mathcal{M} $ at $ x $ is defined as the orthogonal complement of $ T_x\mathcal{M} $ within $ \mathbb{X} $. If given a chart $(\mathcal{U}, \varphi)$ at $x$, a basis of $T_{x} \mathcal{M}$ is given by $\left(\dot{\gamma}_{1}(0),\dot{\gamma}_{2}(0), \cdots, \dot{\gamma}_{n}(0)\right)$, where $ \gamma_{i}(t) := \varphi^{-1}\left(\varphi(x) + t e_{i}\right) $ and $ e_{i} $ denotes the $ i $th canonical vector of $ \mathbb{R}^{n} $. The tangent bundle is defined as $ T\mathcal{M} := \bigcup_{x} T_{x} \mathcal{M} $, representing the collection of all tangent vectors on $ \mathcal{M} $. A map $ V: \mathcal{M} \rightarrow T \mathcal{M} $ is called a vector field on $ \mathcal{M} $ if for every $ x \in \mathcal{M} $, $ V(x) \in T_{x}\mathcal{M} $.

Let $ F: \mathcal{M} \rightarrow \mathbb{X} $ be a smooth mapping. The differential of $ F $ at $ x $, denoted by $ DF(x) $, is a linear map from $ T_{x} \mathcal{M} $ to $ T_{F(x)} \mathbb{X} $. It is defined by $ \left(DF(x)\xi_{x}\right) f := \xi_{x}(f \circ F) $ for any $ \xi_{x} \in T_{x} \mathcal{M} $ and $ f \in \mathfrak{F}_{F(x)}(\mathbb{X}) $. Through the canonical identification $ T_{F(x)}\mathbb{X} \simeq \mathbb{X} $, this differential can be expressed as $ DF(x)\xi_{x} = \sum_{i} (\xi_{x} F^{i}) e_i $, where $ F(x) = \sum_{i} F^{i}(x) e_i $ is the expansion of $ F(x) $ in the basis $ (e_i)_{i=1,2,\cdots,n} $ of $ \mathbb{X} $. If $ \xi = \sum_{i} \xi_{i} \dot{\gamma}_{i}(0) \in T_{x} \mathcal{M} $, then the differential $ D\varphi(x)\xi $ is given by $ \left(\xi_{1}, \xi_{2}, \dots, \xi_{n}\right)^{\top} \in \mathbb{R}^{n} $, meaning that $ D\varphi(x)\xi $ represents the coordinates of $ \xi \in T_{x} \mathcal{M} $ in the chart $ \varphi $. It is known that if $ \mathcal{M} $ is an embedded submanifold of a Euclidean space, the differential $ DF(x) $ reduces to the classical notion of the directional derivative, i.e., $\displaystyle D F(x)\xi_{x}=\lim _{t \rightarrow 0} \frac{F\left(x+t \xi_{x}\right)-F(x)}{t}$.

To distinguish it from Riemannian differential, we use $h'(x)\xi$ to represent the traditional directional derivative for a Euclidean function $h$ in the direction $\xi$ and $\nabla h(x)$ to be the Euclidean gradient of $h$.  We also use these notations for convenience: $	\hat{x}:=\varphi(x)$, $\hat{\xi}:=D \varphi(x)\xi$ and $\hat{f}:=f \circ \varphi^{-1}$. Then for $\xi_{x} \in T_{x} \mathcal{M}$ and $f \in \mathfrak{F}_{x}(\mathcal{M})$, it holds that $	\xi_{x} f=\left\langle\hat{\xi}_{x}, \nabla \hat{f}(\hat{x})\right\rangle_{\mathbb{R}^n}$.

A Riemannian metric $ \langle\cdot, \cdot\rangle_{x} $ is a smoothly varying inner product defined on the tangent space at each point $ x $ of a manifold. A differentiable manifold whose tangent spaces are equipped with such Riemannian metrics is called a Riemannian manifold. Let $ (E_1, E_2, \dots, E_n) $ denote a basis for the vector fields on the manifold, and define $ g_{ij}(x) := \langle E_i(x), E_j(x) \rangle_{x} $. The functions $ g_{ij}(\cdot) $ are then smooth functions on $ \mathcal{M} $. For vector fields $ \xi = \sum_{i} \xi_{i} E_i $ and $ \zeta = \sum_{i} \zeta_{i} E_i $, the inner product at $ x $ is given by
$\langle \xi(x), \zeta(x) \rangle_{x} = \sum_{i,j} g_{ij}(x) \xi^{i}(x) \zeta^{j}(x)$. Let $ G_{\varphi(x)} $ be the matrix-valued function whose $(i, j)$-th entry is $ g_{ij}(x) $. Then, we can obtain
\begin{equation}\label{eq:inner-pro2}
	\langle\xi(x), \zeta(x)\rangle_x=\hat{\xi}(x)^{\top} G_{\varphi(x)} \hat{\zeta}(x).
\end{equation}
Here, $ \hat{\xi}(x) $ and $ \hat{\zeta}(x) $ denote the coordinate representations of the vector fields $ \xi $ and $ \zeta $ at $ x $, respectively. For simplicity, when no confusions arise, we use $ \langle\cdot, \cdot\rangle $ instead of $ \langle\cdot, \cdot\rangle_{x} $. The norm induced by this inner product is denoted $ \|\cdot\| $, with the subscript typically omitted.

Given a function $ f \in \mathfrak{F}_{x}(\mathcal{M}) $, the gradient of $ f $ at the point $ x $, denoted by $ \operatorname{grad} f(x) $, is defined as the unique tangent vector that satisfies
\begin{equation}\label{eq:inner-pro1}
	\langle\operatorname{grad}  f(x), \xi\rangle:=\xi_{x} f=\langle\nabla \hat{f}(\hat{x}), \hat{\xi}\rangle_{\mathbb{R}^n} \quad \forall\, \xi \in T_{x} \mathcal{M}.
\end{equation}
where $ \nabla \hat{f}(\hat{x}) $ is the gradient of the coordinate representation of $ f $ at $ \hat{x} $ in $ \mathbb{R}^n $. The coordinate expression of $ \operatorname{grad} f(x) $ is given by
\begin{equation}\label{eq:co-ex}
	D \varphi(x)\operatorname{grad}  f(x)=G_{\varphi(x)}^{-1} \nabla \hat{f}(\hat{x}).
\end{equation}
where $ G_{\varphi(x)} $ is the matrix whose entries are the components of the Riemannian metric in the given coordinate chart. The Riemannian Hessian of $ f \in \mathfrak{F}_{x}(\mathcal{M}) $ at a point $ x \in \mathcal{M} $ is defined as the symmetric linear map $ \operatorname{Hess} f(x) $ from $ T_{x} \mathcal{M} $ into itself, satisfying $
\operatorname{Hess} f(x)\xi = \nabla_{\xi} \operatorname{grad} f(x)$
for all $ \xi \in T_{x} \mathcal{M} $, where $ \nabla $ is the Riemannian connection on $ \mathcal{M} $. If $ \mathcal{M} $ is embedded in a Euclidean space $ \mathbb{X} $ (see \cite[Section 3.3]{AMS09}), the Riemannian gradient of $ f $ at $ x $ is equal to the orthogonal projection of the Euclidean gradient of $ f $ onto the tangent space $ T_{x} \mathcal{M} $. This relationship is expressed as $
\operatorname{grad} f(x) = \Pi_{x}(\nabla f(x))$, where $ \Pi_{x} $ denotes the projection onto $ T_{x} \mathcal{M} $. Moreover, the Riemannian Hessian of $ f $ at $ x $ in the direction $ \xi $ can be computed as $
\operatorname{Hess} f(x) \xi = \Pi_{x} \left((\operatorname{grad} f)'(x)\xi\right)$,
where $ (\operatorname{grad} f)'(x)\xi $ is the directional derivative of the Riemannian gradient $ \operatorname{grad} f $ at $ x $ in the direction $ \xi $.

The length of a curve $ \gamma: [a, b] \rightarrow \mathcal{M} $ on a Riemannian manifold is defined by $
L(\gamma) = \int_{a}^{b} \sqrt{\langle \dot{\gamma}(t), \dot{\gamma}(t) \rangle} \, \mathrm{d}t$, 
where $ \langle \dot{\gamma}(t), \dot{\gamma}(t) \rangle $ is the Riemannian inner product of the tangent vector $ \dot{\gamma}(t) $ at time $ t $. The Riemannian distance $ d $ between two points $ y, z \in \mathcal{M} $ is then defined as
\begin{equation}\label{eq:def-R-distance}
	d: \mathcal{M} \times \mathcal{M} \rightarrow \mathbb{R}: d(y, z):=\inf _{\Gamma} L(\gamma),	 
\end{equation}
where $ \Gamma $ denotes the set of all curves $ \gamma $ in $ \mathcal{M} $ that connect the points $ y $ and $ z $. Given this distance function, the set $\{y \in \mathcal{M} \mid d(x, y) < \delta\} $ forms a neighborhood of $ x $ with radius $ \delta > 0 $.

With the distance function defined as above, the concept of a Lipschitz function can be extended to a manifold. A function $ f: \mathcal{M} \rightarrow \mathbb{R} $ is said to be Lipschitz of rank $ L > 0 $ on a set $ \mathcal{U} \subseteq \mathcal{M} $ if 
\[
|f(y) - f(z)| \leq L d(y, z) \quad \text{for all } y, z \in \mathcal{U}.
\]
If there exists a neighborhood $ \mathcal{U} $ of $ x \in \mathcal{M} $ such that $ f $ is Lipschitz of rank $ L $ on $ \mathcal{U} $, then $ f $ is said to be Lipschitz of rank $ L $ at $ x $. If $ f $ is Lipschitz of rank $ L $ at every point $ x \in \mathcal{M} $ for some $ L > 0 $, then $ f $ is called locally Lipschitz on $ \mathcal{M} $.

The generalized directional derivative of a locally Lipschitz function $f$ at $x \in \mathcal{M}$ in the direction $v \in T_{x} \mathcal{M}$, is defined in \cite{HP11} as
\begin{equation}\label{eq:direct-div}
	f^{\circ}(x ; v):=\limsup _{y \rightarrow x, t \downarrow 0} \dfrac{f \circ \varphi^{-1}(\varphi(y)+t D \varphi(x)v)-f \circ \varphi^{-1}(\varphi(y))}{t},
\end{equation}
where $(\mathcal{U}, \varphi)$ is a chart containing $x$. {\color{black}
	We should mention that this definition of the generalized directional derivative is independent of the chosen chart around $ x $, as proved in \cite[Lemma 3.3]{MP82}. Therefore, it is straightforward to show that this definition is equivalent to the one provided in \cite{BFO10} by noting that $ \varphi = E^{-1} \operatorname{Exp}_x^{-1} $ also serves as a chart. This chart is commonly known as the normal coordinate chart, which we will discuss in the preliminaries.
}
The Clarke subdifferential of a locally Lipschitz function $f$ at $x \in \mathcal{M}$, denoted by $\partial f(x)$, is defined as
\begin{equation}\label{eq:subdiff}
	\partial f(x)=\left\{\xi \in T_{x} \mathcal{M}\mid\langle\xi, v\rangle \leq f^{\circ}(x ; v) \text { for all } v \in T_{x} \mathcal{M}\right\}.
\end{equation}


The following proposition is the chain rule for the composite function $\theta\circ g$, where $\theta$ is locally Lipschitz and $g$ is continuously differentiable. This extends the result in \cite[Theorem 2.3.10]{C90}.

\begin{proposition}\label{pro:subdiff}
	Suppose that $g:\mathcal{ M }\rightarrow\mathbb{Y}$ is continuously differentiable at $x$ and $\theta:\mathbb{Y}\rightarrow\mathbb{R}$ is locally Lipschitz near $g(x)$. Then $h=\theta\circ g$ is locally Lipschitz near $x$, and one has
	\begin{equation}\label{eq:man-subdiff-include}
	\partial h(x)\subset Dg(x)^{*}\partial\theta(g(x)),
	\end{equation}
{\color{black}	where $Dg(x)^{*}$ is the adjoint operator of $Dg(x)$.} Moreover, the equality holds if $\theta$ is regular (cf. \cite[Definition 2.3.4]{C90}) at $g(x)$. 
\end{proposition}
\proof  The locally Lipschitz of $h$ near $x$ is a direct result of the Lipschitz of $\theta$ and the continuously differentiable of $g$. Suppose that $(\mathcal{ U},\varphi)$ is a chart around $x$. From \cite[Proposition 3.1]{YZS14}, we know that 
\begin{equation}\label{eq:yang-subdiff}
\partial h(x)=(D \varphi(x))^{-1}G_{\varphi(x)}^{-1} \partial\hat{h}(\hat{x}).
\end{equation}
For any $y\in \mathbb{Y} \text{ and } \xi\in T_{x}\mathcal{M}$,  we have
$$
\begin{aligned}
\langle (D \varphi(x))^{-1}G_{\varphi(x)}^{-1}\hat{g}'(\hat{x})^{*}y, \xi\rangle=&\left( D \varphi(x)(D \varphi(x))^{-1}G_{\varphi(x)}^{-1}\hat{g}'(\hat{x})^{*}y\right) ^{\top} G_{\varphi(x)} \left( D \varphi(x) \xi\right) \\
= &\left( \hat{g}'(\hat{x})^{*}y\right) ^{\top}  D \varphi(x) \xi\qquad(\text{by taking } \hat{\xi}=D \varphi(x) \xi) \\
= & \langle  \hat{g}'(\hat{x})^{*}y, \hat{\xi} \rangle	=  \langle \hat{g}'(\hat{x})\hat{\xi},y \rangle=\langle Dg(x)\xi,y \rangle = \langle Dg(x)^*y,\xi\rangle ,
\end{aligned}
$$
which implies that 
\begin{equation}\label{eq:diff-eqv}
Dg(x)^*y=(D \varphi(x))^{-1}G_{\varphi(x)}^{-1}\hat{g}'(\hat{x})^{*}y.
\end{equation}
Since for any $x\in\mathcal{M},\; \partial\theta(g(x))=\partial\theta(\hat{g}(\hat{x}))$, using \cite[Theorem 2.3.10]{C90} we get
\begin{equation}\label{eq:clarke-subdiff}
\partial \hat{h}(\hat{x})\subset \hat{g}'(\hat{x})^{*}\partial\theta(\hat{g}(\hat{x})),
\end{equation}
and the equality holds when $\theta$ is regular.
Combining (\ref{eq:yang-subdiff}), (\ref{eq:diff-eqv}) and (\ref{eq:clarke-subdiff}), we obtain the conclusion.	\qed

{\color{black}\begin{remark}
	By combining \cite[Theorem 3.2]{YZS14} with the fact that the definition of the generalized directional derivative is independent of the choice of chart, we can directly prove Proposition \ref{pro:subdiff}, following a proof strategy similar to that in the Euclidean setting \cite[Theorem 2.3.10]{C90}. For simplicity, we omit the intrinsic proof here.
\end{remark}}

A geodesic on a Riemannian manifold $ \mathcal{M} $ is a curve that locally minimizes the arc length between points. For any $ \xi \in T_{x} \mathcal{M} $, there exists an interval $ \mathcal{I} $ containing zero and a unique geodesic $ \gamma(\cdot; x, \xi): \mathcal{I} \rightarrow \mathcal{M} $ such that $ \gamma(0) = x $ and $ \dot{\gamma}(0) = \xi $. {\color{black} For a sufficiently small neighborhood $\mathcal{U}$ of $0_x$, the mapping
$\operatorname{Exp}_{x}: \mathcal{U} \rightarrow \mathcal{M}, \quad \xi \mapsto \operatorname{Exp}_{x}(\xi) = \gamma(1; x, \xi)$
is called the exponential map at $ x \in \mathcal{M} $. Let $ \mathcal{V} $ be a neighborhood of the origin in $ T_{x} \mathcal{M} $ such that $ \operatorname{Exp}_{x} $ is a diffeomorphism between $ \mathcal{V} $ and $\operatorname{Exp}_{x}(\mathcal{V})$.} Let $ E: \mathbb{R}^{n} \rightarrow T_{x} \mathcal{M} $ be a linear bijection such that the set $ \left\{E(e_{1}), E(e_{2}), \cdots, E(e_{n})\right\} $ forms an orthogonal basis for $ T_{x} \mathcal{M} $, where $ \{e_{1}, e_{2}, \cdots, e_{n}\} $ is the standard basis of $ \mathbb{R}^{n} $. If we define the mapping $ \varphi = E^{-1} \circ \operatorname{Exp}_{x}^{-1} $, then $ (\mathcal{U}, \varphi) $ is known as a Riemannian normal coordinate chart around $ x $.

In a Riemannian normal coordinate chart, the Riemannian metric simplifies, and it is known that $ G_{\varphi(x)} = I_{n} $, where $ I_{n} $ is the identity matrix. As a result, for any function $ f \in \mathfrak{F}_{x}(\mathcal{M}) $ and any tangent vector $ \xi \in T_{x} \mathcal{M} $, the following holds under the normal coordinate chart:
\begin{equation}\label{eq:equi-hessian-inner}
	\langle \xi, \operatorname{Hess}f(x)\xi\rangle=\langle \hat{\xi},\nabla^2_{xx}\hat{f}(\hat{x})\hat{\xi}\rangle
\end{equation}
where $ \hat{\xi} $ is the coordinate representation of $ \xi $, and $ \nabla^2_{xx} \hat{f}(\hat{x}) $ is the Hessian matrix of the coordinate representation of $ f $ at the point $ \hat{x} = \varphi(x) $. This equality expresses the equivalence of the Riemannian Hessian in the manifold and the Euclidean Hessian in the normal coordinates.

In this paper, we are concerned about the continuity and Lipschitz properties of set-valued mappings.
For the set-valued mapping from a Euclidean space $\mathbb{E}$ to a manifold $\mathcal{ M }$, we define the continuity and the (robust) isolated calmness as follows. 
\begin{definition}\label{def:lsc}
	The set-valued mapping $\Psi:\mathbb{E}\rightarrow \mathcal{M}$ is said to be lower semi-continuous at $\bar{p}$ for $\bar{x}$ if for any open neighborhood $\mathcal{V}$ of $\bar{x}$ there exists an open neighborhood $\mathcal{U}$ of $\bar{p}$ such that
	$$
	\emptyset \neq \Psi(p) \cap \mathcal{V} \quad \forall\,p \in \mathcal{U}.
	$$
	The mapping $\Psi$ is said to be upper semi-continuous at $\bar{p}$ if for any open set $\mathcal{O} \supset \Psi(\bar{p})$ there exists an open neighborhood $\mathcal{ U }$ such that for any $p \in \mathcal{U},$ $\Psi(p) \subset \mathcal{O}$. Furthermore, if $\Psi$ is lower semi-continuous at $\bar{p}$ for $\bar{x}$ and is upper semi-continuous at $\bar{p}$, then $\Psi$ is said to be continuous at $(\bar{p}, \bar{x}) \in \operatorname{gph} \Psi$.
\end{definition}

\begin{definition}\label{def:isolated-calm}
	The set-valued mapping $\Psi: \mathbb{E} \rightarrow \mathcal{M}$ is said to be isolated calm at $\bar{p}$ for $\bar{x}$ if there exist a constant $\kappa>0$ and open neighborhoods $\mathcal{U}$ of $\bar{p}$ and $\mathcal{V}$ of $\bar{x}$ such that
	\begin{equation}\label{eq:iso-calm}
		d(x,\bar{x})\leq k\|p-\bar{p}\|\quad\forall\,x\in \Psi(p)\cap\mathcal{V} \quad \mbox{and} \quad p\in \mathcal{U}.
	\end{equation}
\end{definition}
\noindent Moreover, $\Psi$ is said to be robustly isolated calm at $\bar{p}$ for $\bar{x}$ if (\ref{eq:iso-calm}) holds and for each
$p \in \mathcal{U}$, $\Psi(p) \cap \mathcal{V} \neq \emptyset$.


 For a given cone $\mathcal{ C}$, the largest linear space contained in $\mathcal{ C}$ is called the lineality space of $\mathcal{ C}$. A cone $\mathcal{ C}$ is said to be pointed if and only if its lineality space contains only zero. The following $\mathcal{C}^2$-cone reducibility of a closed convex set is taken from \cite[Definition 3.135]{BS13}. 
\begin{definition}\label{def:c2reducible}
	A closed convex set $\mathcal{K}$ is said to be $\mathcal{C}^2$-cone reducible at $y\in \mathcal{K}$, if there exist an open neighborhood $\mathcal{U} \subset \mathbb{Y}$ of $y$, a pointed closed convex cone $\mathcal{C}$ in Euclidean space $\mathbb{Z}$, and a twice continuously differentiable mapping $\Xi: \mathcal{U} \rightarrow \mathbb{Z}$ such that: (i) $\Xi(y)=0 \in \mathbb{Z}$, (ii) the derivative mapping $\Xi^{\prime}(y): \mathbb{Y} \rightarrow \mathbb{Z}$ is onto, and (iii) $\mathcal{K} \cap \mathcal{U}=\{y \in \mathcal{U} \mid \Xi(y) \in \mathcal{C}\} .$ We say that $\mathcal{K}$ is $\mathcal{C }^2$-cone reducible if $\mathcal{K}$ is $\mathcal{C }^2$-cone reducible at every $y \in \mathcal{K}$.
\end{definition}
The class of $ \mathcal{C}^2 $-cone reducible sets is notably rich, including all polyhedral convex sets and many nonpolyhedral sets, such as the second-order cone and positive semidefinite matrices cone \cite{BS13,shapiroSensitivityAnalysisGeneralized2003}. A closed convex function $ \theta: \mathbb{Y} \to \mathbb{R} $ is $ \mathcal{C}^2 $-cone reducible at $ y \in \operatorname{dom}\theta $ if its epigraph $ \operatorname{epi} \theta $ is $ \mathcal{C}^2 $-cone reducible at $ (y, \theta(y)) $. The function $ \theta $ is said to be $ \mathcal{C}^2 $-cone reducible if it is reducible at every $ y \in \operatorname{dom}\theta $.

It follows from \cite[Proposition 3.136]{BS13} that if $\mathcal{D}$ is a $\mathcal{C}^2$-cone reducible convex set, then $\mathcal{T}_{\mathcal{D}}^{i,2}(y,h)=\mathcal{T}_{\mathcal{D}}^{2}(y, h)$. In this case, $\mathcal{T}_{\mathcal{D}}^{2}(y, h)$ will be simply called the second order tangent set to $\mathcal{D}$ at $y \in \mathcal{D}$ in the direction $h \in \mathbb{Y}$.

Let $\mathcal{D}$ be a closed set in $\mathbb{Y}$. The radial cone at a point $y \in \mathcal{D}$ is defined by
\begin{equation}\label{eq:def-radialcone}
	\mathcal{R}_{\mathcal{D}}(y):=\left\{d \in \mathbb{Y} \mid \exists t^{*}>0 \text { such that } y+t d \in \mathcal{D} \text { for any } t \in\left[0, t^{*}\right]\right\},
\end{equation}
and the tangent cone at $y$ is defined as
$$
\mathcal{T}_{\mathcal{D}}(y)=\left\{d \in \mathbb{Y} \mid \exists t^{k} \downarrow 0,\; \operatorname{dist}\left(y+t^{k} d, \mathcal{D}\right)= o\left(t^{k}\right)\right\},
$$
where $\operatorname{dist}(u, \mathcal{D}) = \inf \{ \| u - p \| \mid p \in \mathcal{D} \}$ denotes the Euclidean distance from $u$ to $\mathcal{D}$.  The normal cone to $\mathcal{D}$ at $y$ is defined as $
\mathcal{N}_{\mathcal{D}}(y) = \left( \mathcal{T}_{\mathcal{D}}(y) \right)^{\circ}$, 
which is the polar cone of $\mathcal{T}_{\mathcal{D}}(y)$. The inner and outer second-order tangent sets to the closed set $\mathcal{D}$ in the direction $h \in \mathbb{Y}$ are defined respectively by
$$
\mathcal{T}_{\mathcal{D}}^{i, 2}(y, h):=\left\{w \in \mathbb{Y} \mid \operatorname{dist}\left(y+t h+\dfrac{1}{2} t^{2} w, \mathcal{D}\right)=o\left(t^{2}\right), t \geq 0\right\},
$$
and
$$
\mathcal{T}_{\mathcal{D}}^{2}(y, h):=\left\{w \in \mathbb{Y} \mid \exists t_{k} \downarrow 0 \text { such that } \operatorname{dist}\left(y+t_{k} h+\dfrac{1}{2} t_{k}^{2} w, \mathcal{D}\right)=o\left(t_{k}^{2}\right)\right\}.
$$
According to \cite[Proposition 3.136]{BS13}, if $\mathcal{D}$ is a $\mathcal{C}^2$-cone reducible convex set, then the inner and outer second-order tangent sets coincide: $\mathcal{T}_{\mathcal{D}}^{i, 2}(y, h) = \mathcal{T}_{\mathcal{D}}^{2}(y, h)$. In such cases, $\mathcal{T}_{\mathcal{D}}^{2}(y, h)$ is simply referred to as the second-order tangent set to $\mathcal{D}$ at $y$ in the direction $h$.

For a given function $\theta: \mathbb{Y} \rightarrow(-\infty,+\infty]$, the lower and upper directional epiderivatives (cf. e.g., \cite[(2.68) and (2.69)]{BS13}) of $\theta$ at $y \in \operatorname{dom} \theta$ in the direction $h \in \mathbb{Y}$ are defined as
$$
\theta_{-}^{\downarrow}(y ; h) :=\liminf _{t \downarrow 0 \atop h^{\prime} \to h} \frac{\theta\left(y+t h^{\prime}\right)-\theta(y)}{t} \quad {\rm and} \quad \theta_{+}^{\downarrow}(y ; h) :=\sup _{\left\{t_{n}\right\} \in \Sigma}\left(\liminf _{n \to \infty \atop h^{\prime} \to h} \frac{\theta\left(y+t_{n} h^{\prime}\right)-\theta(y)}{t_{n}}\right),
$$
respectively, where $\Sigma$ is the set of all positive sequences $\left\{t_{n}\right\}$ converging to zero. Since $\theta$ is a proper closed convex function, it follows that $\theta_{-}^{\downarrow}(y ; \cdot)=\theta_{+}^{\downarrow}(y ; \cdot)$. Moreover, it is known \cite[Proposition 2.58]{BS13} that for any $y\in \operatorname{dom} \theta$,
\begin{equation}\label{eq:tangent-epideri}
	\mathcal{T}_{\operatorname{epi} \theta}(y, \theta(y))=\operatorname{epi} \theta_{-}^{\downarrow}(y ; \cdot).
\end{equation}
If $\theta_{-}^{\downarrow}(y; h)$ is finite for $x \in \operatorname{dom} \theta$ and $h \in \mathbb{Y}$, the lower second order epiderivatives \cite[(2.76)]{BS13} for $w\in \mathbb{Y}$ is defined as:
$$
\theta_{-}^{\downarrow \downarrow}(y ; h, w):=\liminf _{t \downarrow 0 \atop w^{\prime} \to w} \frac{\theta\left(y+t h+\frac{1}{2} t^{2} w^{\prime}\right)-\theta(y)-t \theta_{-}^{\downarrow}(y ; h)}{\frac{1}{2} t^{2}}.
$$
Finally, it follows from \cite[Proposition 3.41]{BS13} that
\begin{equation}\label{eq:2d-tan}
	\mathcal{T}_{\operatorname{epi} \theta}^{2}\left((y, \theta(y)) ;h, \theta_{-}^{\downarrow}(y;h)\right)=\operatorname{epi} \theta_{-}^{\downarrow \downarrow}(y ; h, \cdot).
\end{equation}

\section{The robust isolated calmness of {\color{black}the} KKT solution mapping}\label{sec:ric}
In this section, we study the robust isolated calmness of the KKT solution mapping for problem~(\ref{eq:prime}). We will introduce the Robinson constraint qualification, the strict Robinson constraint qualification, and the second-order optimality conditions for Riemannian manifolds in the following subsection. We will then define the perturbed problem using the normal coordinate chart and characterize the robust isolated calmness property for (\ref{eq:prime}).

\subsection{The constraint qualifications and second order optimality conditions}\label{subsec:cq}
Consider the generalized form of the optimization problem on manifold as follows:
\begin{equation}\label{eq:man-cone}
	\begin{array}{ll}
		\min & f(x) \\
		\text { s.t. } 
		& g(x) \in \mathcal{K},\\
		& x\in \mathcal{M},
	\end{array}
\end{equation}
where $f: \mathcal{M} \rightarrow \mathbb{R}$ and $ g: \mathcal{M} \rightarrow \mathbb{Y}$ are twice continuously differentiable, and $\mathcal{K} \subset \mathbb{Y}$ is a nonempty closed convex set  in $\mathbb{Y}$. When $g(x)=(g_1(x),g_2(x))$ and $\mathcal{K}=\operatorname{epi}\theta\times\mathcal{Q}$, problem (\ref{eq:man-cone}) is equivalent to the nonsmooth problem (\ref{eq:prime}). In addition, we also assume $\mathcal{K}$ is $C^{2}$-cone reducible (Definition \ref{def:c2reducible}). Using the normal coordinate chart, we can locally transform (\ref{eq:man-cone}) into the traditional minimizing problem in Euclidean space. For $x\in \mathcal{F}:=\{x \in \mathcal{M} \mid g(x) \in \mathcal{K}\}$, suppose $(\mathcal{U},\varphi)$ is the normal coordinate chart around $x$. Then, we locally obtain an equivalent problem in $\mathbb{R}^n$ 
\begin{equation}\label{eq:rn-cone}
	\begin{array}{ll}
		\min & \hat{f}(\hat{x}) \\
		\text { s.t. } & \hat{g}(\hat{x})\in \mathcal{K},\\
		&\hat{x}\in \varphi(\mathcal{U})\in \mathbb{R}^n,
	\end{array} 
\end{equation}
where $\hat{f}=f\circ \varphi^{-1}$, $\hat{g}=g\circ \varphi^{-1}$ and $\hat{x}=\varphi(x)=(\varphi_1(x),\varphi_2(x),\cdots,\varphi_n(x))$. Suppose that $\hat{x}$ is a feasible solution to problem (\ref{eq:rn-cone}). The critical cone $\widehat{\mathcal{C}}(\hat{x})$ of (\ref{eq:rn-cone}) at $\hat{x}$ is defined by
\begin{equation}\label{eq:rn-criti}
	\widehat{\mathcal{C}}(\hat{x}):=\left\{\hat{d} \in \mathbb{R}^n \mid \hat{g}'(\hat{x}) \hat{d} \in \mathcal{T}_{\mathcal{K}}(\hat{g}(\hat{x})), \hat{f}'(\hat{x}) \hat{d} \leq 0\right\}.
\end{equation}
For the corresponding $x\in \mathcal{M}$ for problem (\ref{eq:man-cone}), we can define the critical cone $\mathcal{C}(x)$ of (\ref{eq:man-cone}) at $x$ as
\begin{equation}\label{eq:man-cri}
	\mathcal{C}(x):=\left\{d \in T_x\mathcal{M} \mid Dg(x) d \in \mathcal{T}_{\mathcal{K}}(g(x)), Df(x) d \leq 0\right\}.
\end{equation}
It is easy to see that $\widehat{\mathcal{C}}(\hat{x})=D\varphi(x)\mathcal{C}(x)$.

Let $L: \mathcal{M} \times \mathbb{Y} \rightarrow \mathbb{R}$ be the Lagrangian function of problem (\ref{eq:man-cone}) defined by
\begin{equation}\label{eq:man-lagran}
	L(x ; y):=f(x)+\langle y, g(x)\rangle, \quad(x, y) \in \mathcal{M} \times \mathbb{Y}.
\end{equation}
and $\widehat{L}:\mathbb{R}^n\times\mathbb{Y}\rightarrow \mathbb{R}$ be the Lagrangian function of problem (\ref{eq:rn-cone}) defined by
\begin{equation}\label{eq:rn-lagran}
	\widehat{L}(\hat{x}; y):=\hat{f}(\hat{x})+\langle y, \hat{g}(\hat{x})\rangle, \quad(\hat{x}, y) \in \mathbb{R}^n  \times \mathbb{Y}.
\end{equation}
We say that $x\in\mathcal{M}$ is a stationary point of (\ref{eq:man-cone}) and $y \in \mathbb{Y}$ is a Lagrange multiplier  at  $x$, if $(x,y)$ satisfies the Karush-Kuhn-Tucker (KKT) condition:
\begin{equation}\label{eq:man-kkt}
	\left\{\begin{array} { l } 
		{\text{grad}_xL(x;y) =0 , } \\
		y \in \mathcal{N}_{\mathcal{K}}(g ( x )),
	\end{array}\right.
\end{equation}
where $\mathcal{N}_{\mathcal{K}}(g(x))$ is the normal cone to $\mathcal{K}$ at $g(x) \in \mathbb{Y}$. We denote by $M(x)$ the set of all Lagrange multipliers at $x$. We also use $ \widehat{M}(\hat{x})$ to denote the set of multipliers at a stationary point $\hat{x}$ for problem (\ref{eq:rn-cone}).

It is well known (cf. e.g., \cite{BS13}) that the Robinson constraint qualification (RCQ) for problem (\ref{eq:rn-cone}) holds at a feasible point $\hat{x}$ if
\begin{equation}\label{eq:rn-rcq}
	\hat{g}'(\hat{x}) \mathbb{R}^n+\mathcal{T}_{\mathcal{K}}(\hat{g}(\hat{x}))=\mathbb{Y},
\end{equation}
{\color{black}
	where for any $ d \in \hat{g}'(\hat{x}) \mathbb{R}^n $, there exists some $ v \in \mathbb{R}^n $ such that $ d = \hat{g}'(\hat{x})[v] $. Moreover, the sum in \eqref{eq:man-kkt} is interpreted as the sum of sets, meaning that  $$
	\hat{g}'(\hat{x}) \mathbb{R}^n + \mathcal{T}_{\mathcal{K}}(\hat{g}(\hat{x})) = \{ d_1 + d_2 \mid d_1 \in \hat{g}'(\hat{x}) \mathbb{R}^n, \ d_2 \in \mathcal{T}_{\mathcal{K}}(\hat{g}(\hat{x})) \}.$$
	}The strict Robinson constraint qualification (SRCQ) for problem (\ref{eq:rn-cone}) holds with respect to $y\in \widehat{M}(\hat{x})$ if for the stationary point $\hat{x}$ if
\begin{equation}\label{eq:rn-srcq}
	\hat{g}'(\hat{x}) \mathbb{R}^n+\mathcal{T}_{\mathcal{K}}(\hat{g}(\hat{x})) \cap y^{\perp}=\mathbb{Y},
\end{equation}
and the constraint non-degeneracy is said to hold at $\hat{x}$ if
$$
\hat{g}'(\hat{x}) \mathbb{R}^n + \operatorname{lin}\left(\mathcal{T}_{\mathcal{K}}(\hat{g}(\hat{x}))\right)=\mathbb{Y},
$$
where $\operatorname{lin}\left(\mathcal{T}_{\mathcal{K}}(g(\hat{x}))\right)$ is the lineality space of $\mathcal{T}_{\mathcal{K}}(\hat{g}(\hat{x}))$. Therefore, we are able to define the Robinson constraint qualification, the strict Robinson constraint qualification, and the constraint non-degeneracy on Riemannian manifolds. These definitions can lead us to the existence or uniqueness of Lagrange multipliers for problem (\ref{eq:man-cone}).

\begin{definition}\label{def:rcq+srcq}
	For the manifold optimization problem \eqref{eq:man-cone}, we say that the manifold Robinson constraint qualification (M-RCQ) holds at $x\in \mathcal{ F }$ if 
	\begin{equation}\label{eq:man-rcq}
		Dg(x) T_{x} \mathcal{M}+\mathcal{T}_{\mathcal{K}}(g(x))=\mathbb{Y},
	\end{equation}
	the manifold strict Robinson constraint qualification (M-SRCQ) holds at a stationary point $x$ with respect to $y\in M(x)$ if 
	\begin{equation}\label{eq:man-srcq}
		Dg(x) T_{x} \mathcal{M}+\mathcal{T}_{\mathcal{K}}(g(x)) \cap y^{\perp}=\mathbb{Y},
	\end{equation}
	and the manifold constraint non-degeneracy  holds at $x\in \mathcal{ F }$ if 
	\begin{equation}\label{eq:con-nondegen}
		Dg(x) T_{x} \mathcal{M}+\operatorname{lin}\left(\mathcal{T}_{\mathcal{K}}(g(x))\right)=\mathbb{Y}.
	\end{equation}
\end{definition}

{\color{black}
In Euclidean spaces, the RCQ ensures the existence of Lagrange multipliers, while the SRCQ, a stronger condition, guarantees the uniqueness of these multipliers. Additionally, constraint non-degeneracy ensures that the gradients of the active constraints at the optimal solution are linearly independent. For Euclidean NLPs, it is known that the RCQ and constraint non-degeneracy conditions correspond to the MFCQ and LICQ, respectively \cite[page 71 and Example 4.78]{BS13}. In the manifold setting, our definitions of M-RCQ and manifold constraint non-degeneracy are also equivalent to the manifold MFCQ and manifold LICQ conditions proposed in \cite{BH19,YZS14} when $ \mathcal{K} $ is polyhedral. In the following discussions, we will show that M-RCQ, M-SRCQ, and manifold constraint non-degeneracy play crucial roles in ensuring the existence and uniqueness of the KKT solution.
	}

By using the normal coordinate chart, we establish the following theorem on the existence and boundedness of multipliers for problem (\ref{eq:man-cone}).

\begin{theorem}\label{thm:man-rcq}
	Suppose that $x^*$ is a locally optimal solution of (\ref{eq:man-cone}). Then, $M(x^*)$ is a nonempty, convex, bounded, and compact subset of $\mathbb{Y}$ if and only if the M-RCQ (\ref{eq:man-rcq}) holds at $x^*$.
\end{theorem}
\proof Recall that $\widehat{M}(\hat{x})$ is the set of all Lagrange multipliers at $\hat{x}$ for problem (\ref{eq:rn-cone}). Then, for any $y\in\widehat{M}(\hat{x})$, $y$ satisfies the KKT condition 
\begin{equation}\label{eq:rn-kkt}
	\left\{\begin{array} { l } 
		\nabla_{\hat{x}} \widehat{L}(\hat{x};y) =0 , \\
		y \in \mathcal{N}_{\mathcal{K}}(\hat{g} ( \hat{x} )).
	\end{array}\right.
\end{equation}
It then follows from (\ref{eq:co-ex}) that $\operatorname{grad} _{x}L(x,y)=(D\varphi(x))^{-1}G_{\varphi(x)}^{-1}\nabla_{ \hat{x} }\widehat{L}(\hat{x},y)=0$, which implies that $y\in M(x).$ Thus, $ \widehat{M}(\hat{x})\subseteq M(x)$. The inverse relation $M(x)\subseteq \widehat{M}(\hat{x})$ can be obtained similarly. Hence, we have $M(x)= \widehat{M}(\hat{x})$.

Suppose that $x^*$ is a locally optimal solution of (\ref{eq:man-cone}), then $\hat{x}^*$ is a local solution of the equivalent problem (\ref{eq:rn-cone}). Thus, for any $d\in T_{x^*}\mathcal{ M }$, by taking $\hat{d}=D\varphi(x^*)d$, we obtain that $Dg(x^*)d=\hat{g}'(\hat{x}^*)\hat{d}$. Conversely, for any given $\hat{d}\in \mathbb{R}^n$, $d=(D\varphi(x^*))^{-1}\hat{d}$ is a tangent vector at $x^*$, which ensures that $Dg(x^*)d=\hat{g}'(\hat{x}^*)\hat{d}$. Therefore, the M-RCQ (\ref{eq:man-rcq}) holds at $x^*$ if and only if RCQ (\ref{eq:rn-rcq}) holds at $\hat{x}^*$.
The conclusion is then obtained by the well-known fact \cite{zoweRegularityStabilityMathematical1979} that $\widehat{M}(\hat{x}^*)$ is nonempty, convex, bounded and compact if and only if RCQ (\ref{eq:rn-rcq}) holds at $\hat{x}^*$.\qed

\begin{remark}\label{rm:rcq-necessary-localmin}
	By Theorem \ref{thm:man-rcq}, the M-RCQ (\ref{eq:man-rcq}) is a necessary condition for a given stationary point $x^*$ to be the local optimal solution of (\ref{eq:man-cone}). This result extends the first-order necessary condition for 
$\mathcal{K}$ being polyhedral, as proposed in \cite[Theorem 4.1]{BH19}, to a more general context of manifold-constrained optimization.
\end{remark}

The statement in Theorem \ref{thm:man-srcq} can also be derived using the normal coordinate chart. This is shown by proving that M-SRCQ (\ref{eq:man-srcq}) at $ x^* $ with respect to $ y^* $ is equivalent to SRCQ (\ref{eq:rn-srcq}) at $ \hat{x}^* $ for $ y^* $. However, we provide a direct proof here, following the approach in \cite[Proposition 4.47]{BS13}.

\begin{proposition}[\textbf{Uniqueness of Lagrangian multipliers}]\label{pro:uni-multi}
	Suppose that $x^*$ is a feasible solution to (\ref{eq:man-cone}) and $M(x^*) \neq \emptyset .$ Let $y^* \in M(x^*) .$ Then, the corresponding multiplier set $M(x^*)$ is a singleton if and only if
	\begin{equation}\label{eq:man-m(x)sin}
		\left[Dg(x^*) T_{x^*}\mathcal{M}\right]^{\perp} \cap \mathcal{R}_{\mathcal{N}_{\mathcal{K}}(g(x^*))}\left(y^*\right)=\{0\} .
	\end{equation}
\end{proposition}
\proof For any point $ y \in M(x^*) $, let $ \lambda = y - y^* $. We know that $ Dg(x^*)^* \lambda = 0 $ and $y^*+\lambda \in \mathcal{N}_{\mathcal{K}}(g(x^*))$. Since $ y^* \in \mathcal{N}_{\mathcal{K}}(g(x^*)) $, the definition of the radial cone \eqref{eq:def-radialcone} implies that $ \lambda \in \mathcal{R}_{\mathcal{N}_{\mathcal{K}}(g(x^*))}(y^*) $. The condition $ Dg(x^*)^* \lambda = 0 $ indicates that for any $ \xi \in T_{x^*} \mathcal{M} $, we have $ \langle Dg(x^*) \xi, \lambda \rangle = 0 $. Consequently, $ \lambda $ is orthogonal to the image of $ Dg(x^*) $ in $ T_{x^*} \mathcal{M} $, meaning $ \lambda \in \left[ Dg(x^*) T_{x^*} \mathcal{M} \right]^{\perp} $. Moreover, for any $ \tilde{\lambda} \in \left[ Dg(x^*) T_{x^*} \mathcal{M} \right]^{\perp} $, it follows that $ 0 = \langle Dg(x^*) \xi, \tilde{\lambda} \rangle = \langle \xi, Dg(x^*)^* \tilde{\lambda} \rangle $. Thus, $ Dg(x^*)^* \tilde{\lambda} = 0 $.Therefore, $ M(x^*) $ is not a singleton if and only if there exists a nonzero $ \lambda $ such that $ \lambda $ lies in both $ \left[ Dg(x^*) T_{x^*} \mathcal{M} \right]^{\perp} $ and $ \mathcal{R}_{\mathcal{N}_{\mathcal{K}}(g(x^*))}(y^*) $. This concludes the proof. \qed

By Proposition \ref{pro:uni-multi}, we obtain the following result on the uniqueness of Lagrangian multipliers of (\ref{eq:man-cone}), immediately.
\begin{theorem}\label{thm:man-srcq}
	Suppose that $x^*$ is a feasible solution to (\ref{eq:man-cone}) and $M(x^*) \neq \emptyset .$ Let $y^* \in M(x^*)$, then $M(x^*)$ is a singleton if M-SRCQ (\ref{eq:man-srcq}) holds at $x^*$ with respect to $y^*$.
\end{theorem}
\proof Assume that the M-SRCQ condition (\ref{eq:man-srcq}) holds at $(x^*, y^*)$. 
Let $S=D g(x^*) T_{x^*} \mathcal{M}+\mathcal{T}_{\mathcal{K}}(g(x^*))$. Then, we have 
$$
S^{\circ}=\left(Dg(x^*) T_{x^*}\mathcal{M}\right)^{\circ} \cap\left(\mathcal{T}_{\mathcal{K}}(g(x^*)) \cap {y^*}^{\perp}\right)^{\circ} .
$$
Since $\left(Dg(x^*) T_{x^*}\mathcal{M}\right)^{\circ}=\left(Dg(x^*) T_{x^*}\mathcal{M}\right)^{\perp}$ and
$$
\left(\mathcal{T}_{\mathcal{K}}(g(x^*)) \cap {y^*}^{\perp}\right)^{\circ}=\operatorname{cl}\left(\mathcal{N}_{\mathcal{K}}(g(x^*))+\operatorname{span}\left(y^*\right)\right)=\mathcal{T}_{\mathcal{N}_{\mathcal{K}}(g(x^*))}\left(y^*\right),
$$
it follows that 
$$S^{\circ}=\left(Dg(x^*) T_{x^*}\mathcal{M}\right)^{\perp} \cap \mathcal{T}_{\mathcal{N}_{\mathcal{K}}(g(x^*))}\left(y^*\right).$$
By noting that $\mathcal{R}_{\mathcal{N}_{\mathcal{K}}(g(x^*))}\left(y^*\right) \subset \mathcal{T}_{\mathcal{N}_{\mathcal{K}}(g(x^*))}\left(y^*\right)$, we know that (\ref{eq:man-m(x)sin}) is satisfied, which implies that $M(x^*)=\{y^*\}$. This completes the proof.\qed

\begin{remark}\label{rm:con-nonde-singleton}
It is straightforward to observe that the manifold constraint non-degeneracy condition is stronger than the M-SRCQ condition. This is because $\operatorname{lin}(\mathcal{T}_{\mathcal{K}}(g(x^*))) \subset \mathcal{T}_{\mathcal{K}}(g(x^*)) \cap {y^*}^{\perp}$, as shown in the proof of \cite[Proposition 4.73]{BS13}. Consequently, constraint non-degeneracy (\ref{eq:con-nondegen}) also guarantees the uniqueness of the multiplier.
\end{remark}

\begin{remark}\label{no-chart-relying-rcq}
It is worth noting that the equivalence between the M-RCQ or M-SRCQ for problem (\ref{eq:man-cone}) and the Euclidean RCQ or SRCQ for (\ref{eq:rn-cone}) holds in any chart around the given point. For simplicity, we will use the normal coordinate chart for this discussion.
\end{remark}

\begin{remark}\label{rm:em-srcq-equi}
	When $\mathcal{ M}$ reduce to a special kind of embedded manifold which is written as $\mathcal{ M}=\left\lbrace x\in\mathbb{X}\mid h(x)=0 \right\rbrace $, where $h:\mathbb{X}\rightarrow\mathbb{R}^{m}$ is smooth and $h'(x)$ has full row rank for all $x\in\mathbb{X}$, then by \cite[Section 3.5.7]{AMS09} the tangent space at $x\in\mathcal{M}$ is  $T_x\mathcal{ M}=\operatorname{Ker}(h'(x))$. If we regard (\ref{eq:man-cone}) as a constrained Euclidean optimization, the Euclidean RCQ condition at a feasible point $x$ for this problem is given by
	\begin{equation}\label{eq:em-rcq}
		\begin{bmatrix}
			g'(x)\\
			h'(x)
		\end{bmatrix}\mathbb{X}+\begin{bmatrix}
			\mathcal{ T}_{\mathcal{K}}\left( g(x)\right)\\
			\mathcal{ T}_{\left\lbrace 0 \right\rbrace^{m} }\left( h(x)\right) \\
		\end{bmatrix}=\begin{bmatrix}
			\mathbb{Y}\\
			\mathbb{R}^{m}
		\end{bmatrix},
	\end{equation}
	the SRCQ condition at $x$ with respect to the corresponding multiplier $(y,z)$ is 
	\begin{equation}\label{eq:em-srcq}
		\begin{bmatrix}
			g'(x)\\
			h'(x)
		\end{bmatrix}\mathbb{X}+\begin{bmatrix}
			\mathcal{ T}_{\mathcal{K}}\left( g(x)\right)\cap y^{\perp}\\
			\mathcal{ T}_{\left\lbrace 0 \right\rbrace^{m} }\left( h(x)\right)\cap z^{\perp} \\
		\end{bmatrix}=\begin{bmatrix}
			\mathbb{Y}\\
			\mathbb{R}^{m}
		\end{bmatrix},
	\end{equation}
	and the constraint non-degeneracy at $x$ is 
	\begin{equation}\label{eq:em-consnondegeneracy}
		\begin{bmatrix}
			g'(x)\\
			h'(x)
		\end{bmatrix}\mathbb{X}+\begin{bmatrix}
			\operatorname{lin}\left(\mathcal{T}_{\mathcal{K}}(g(x))\right)\\
			\operatorname{lin}\left(	\mathcal{ T}_{\left\lbrace 0 \right\rbrace^{m} }\left( h(x)\right)\right) \\
		\end{bmatrix}=\begin{bmatrix}
			\mathbb{Y}\\
			\mathbb{R}^{m}
		\end{bmatrix}.
	\end{equation}
	We will show that the RCQ condition (\ref{eq:em-rcq})  (SRCQ (\ref{eq:em-srcq}), constraint non-degeneracy (\ref{eq:em-consnondegeneracy})) is equivalent to the M-RCQ  (M-SRCQ, manifold constraint non-degeneracy) condition.

	It is obvious that $\mathcal{ T}_{\left\lbrace 0 \right\rbrace^{m} }\left( h(x)\right)=\{0\}$ for any feasible $x$. If M-RCQ  (\ref{eq:man-rcq}) holds at $x$, then RCQ  (\ref{eq:em-rcq}) is satisfied by $Dg(x)T_x{\mathcal{M}}=g'(x)T_x\mathcal{M}\subseteq g'(x)\mathbb{X}$ 
	and the fact that $h'(x)$ is of full row rank. Conversely, if  RCQ  (\ref{eq:em-rcq}) holds at $x$, for any $(d_1,d_2)\in \mathbb{Y}\times\mathbb{R}^m$, there exist $\xi \in \mathbb{X}$ and $\eta\in \mathcal{T}_{\mathcal{K}}(g(x))$, such that $g'(x)\xi+\eta=d_1$ and $h'(x)\xi=d_2$. Specially taking $d_2=0$, we have $\xi \in \operatorname{Ker}\left(h'(x)\right)=T_x\mathcal{M},$ which implies that M-RCQ (\ref{eq:man-rcq}) is fulfilled. The equivalence between the constraint non-degeneracy conditions can be similarly obtained. For the SRCQ conditions, we only need to show that at any stationary point $x$, the unique multiplier to the constraint of (\ref{eq:man-cone}) is the same as the multiplier to the Euclidean form problem.	The KKT condition of the constrained Euclidean problem can be written as 
	\begin{equation}\label{eq:em-kkt}
		\left\{\begin{array} { l } 
			{\tilde{L}(x,y,z)=\nabla f(x)+g'(x)^{*}y+h'(x)^{*}z =0 , } \\
			y \in \mathcal{N}_{\mathcal{K}}(g ( x )),\\
			h(x)=0.
		\end{array}\right.
	\end{equation}
	The first equality in (\ref{eq:em-kkt}) implies that $\nabla f(x)+g'(x)^{*}y\in \operatorname{Range}(h'(x)^{*})=N_x\mathcal{M}$, or equivalently, 
	$\operatorname{grad}  f(x)+Dg(x)^{*}y = \Pi_x(\nabla f(x)+g'(x)^{*}y)=0$. Thus the SRCQ condition (\ref{eq:em-srcq}) holds at $(x,y,z)$ if and only if  the M-SRCQ condition (\ref{eq:man-srcq}) holds at $(x,y)$.
	
	It is worth mentioning that when $h'(x)$ is not onto, the RCQ condition (\ref{eq:em-rcq}) cannot be satisfied, since $h'(x)\mathbb{X}\neq \mathbb{R }^{m}$ and $\mathcal{T}_{\{0\}^m}(h(x))=\{0\}$. Therefore, the SRCQ condition (\ref{eq:em-srcq}) and the constraint non-degeneracy (\ref{eq:em-consnondegeneracy}) can also never hold. 
	

\end{remark}

For the $\mathcal{C}^2$-cone reducible set $\mathcal{ K}$, we will introduce the second order necessary and sufficient conditions to problem (\ref{eq:man-cone}). Firstly, we define the quadratic growth condition.

\begin{definition}\label{def:man-2grow}
	We say that the quadratic growth condition holds at $x^*\in\mathcal{ M }$, if there exist a neighborhood $\mathcal{ N }$ of $x^*$ and a positive constant $c>0$ such that for all $x\in\mathcal{ F }\cap\mathcal{ N}$, the following inequality holds: 
	$$
	f(x)\geq f(x^*)+cd^2(x,x^*).
	$$ 
\end{definition}

{\color{black}The quadratic growth condition is important in optimization because it ensures that the objective function grows quadratically as the solution moves away from the optimal point. This condition guarantees that the optimization problem has a unique solution within a neighborhood of the optimum. Moreover, its connection with the second order sufficient condition is shown in the following Theorem \ref{thm:man-sosc}.}

\begin{theorem}\label{thm:man-sosc}
	Suppose that $x^*$ is a locally optimal solution to problem (\ref{eq:man-cone}) and M-RCQ (\ref{eq:man-rcq}) holds at $x^*$. Then the following manifold second order necessary  condition (M-SONC) holds:
	\begin{equation}\label{eq:man-sonc}
		\sup _{y \in M(x^*)}\left\{\left\langle \xi, \operatorname{Hess}_{x} L(x^* ; y) \xi\right\rangle-\sigma\left(y, \mathcal{T}_{\mathcal{K}}^{2}\left(g(x^*), Dg(x^*) \xi\right)\right)\right\} \geq 0 \quad \forall\,\xi \in \mathcal{C}(x^*),
	\end{equation}
	where for any $y \in \mathbb{Y}$, $\operatorname{Hess}_{x} L(x^* ; y)$ is the Hessian of $L(\cdot ; y)$ at $x^*$ and $\sigma(y,\mathcal{D})$ is the support function of set $\mathcal{D}$ at $y$. Conversely, suppose $x^*$ is a stationary point of problem (\ref{eq:man-cone}) and M-RCQ (\ref{eq:man-rcq}) holds at $x^*$. Then the following manifold second order sufficient condition (M-SOSC)
	\begin{equation}\label{eq:man-sosc}
		\sup _{y \in M(x^*)}\left\{\left\langle \xi, \operatorname{Hess}_{x} L(x^* ; y) \xi\right\rangle-\sigma\left(y, \mathcal{T}_{\mathcal{K}}^{2}\left(g(x^*), Dg(x^*) \xi\right)\right)\right\}>0 \quad \forall\,\xi \in \mathcal{C}(x^*) \backslash\{0\}
	\end{equation}
	is necessary and sufficient for the quadratic growth condition at the point $x^*$ for problem (\ref{eq:man-cone}).
\end{theorem}
\proof Let $(\mathcal{U},\varphi)$ be the normal coordinate chart at $x^*$.  Note that M-RCQ (\ref{eq:man-rcq}) holds at $x^*$ indicates that RCQ (\ref{eq:rn-rcq}) holds at $\hat{x}^*$. It is well known \cite[Theorem 5.2]{B99} in Euclidean space that the traditional second order necessary condition holds for problem (\ref{eq:rn-cone}) at the locally optimal solution $\hat{x}^*=\varphi(x^*)$ if RCQ (\ref{eq:rn-rcq}) holds. Thus, by (\ref{eq:equi-hessian-inner}), we have
$$
\left\langle \xi, \operatorname{Hess}_{x} L(x^* ; y) \xi\right\rangle=\left\langle \hat{\xi}, \nabla_{x x}^2 \widehat{L}(\hat{x}^* ; y) \hat{\xi}\right\rangle \quad\forall\,\xi\in \mathcal{C}(x^*)\quad \mbox{and}\quad \hat{\xi}=D\varphi(x^*)\xi,
$$
which implies that the inequality (\ref{eq:man-sonc}) holds if $x^*$ is a local solution and M-RCQ  (\ref{eq:man-rcq}) holds. 

Conversely, if $x^*$ is a stationary point of problem (\ref{eq:man-cone}), then $\varphi(x^*)$ is a stationary point of problem (\ref{eq:rn-cone}). Therefore, the SOSC condition (\ref{eq:rn-sosc}) for problem (\ref{eq:rn-cone}), which is stated as below, is the necessary and sufficient condition for the quadratic growth condition at $\varphi(x^*)$ for problem~(\ref{eq:rn-cone}) 
\begin{equation}\label{eq:rn-sosc}
	\sup _{y \in \widehat{M}(\hat{x}^*)}\left\{\left\langle \hat{\xi}, \nabla^2_{\hat{x} \hat{x}} \widehat{L}(\hat{x}^*; y) \hat{\xi}\right\rangle-\sigma\left(y, \mathcal{T}_{\mathcal{K}}^{2}\left(\hat{g}(\hat{x}^*), \hat{g}'(\hat{x}^*) \hat{\xi}\right)\right)\right\}>0 \quad \forall\,\hat{\xi} \in \widehat{\mathcal{C}}(\hat{x}^*) \backslash\{0\}.
\end{equation}
The quadratic growth condition for (\ref{eq:rn-cone}) at $\hat{x}^*$ is that there exist a neighborhood $\widehat{\mathcal{ N}}$ of $\hat{x}^*$ and a positive constant $c>0$, such that for all $\hat{x}\in \mathcal{F }\cap \widehat{\mathcal{ N}}$, the following inequality holds 
\begin{equation}\label{eq:rn-2grow}
	\hat{f}(\hat{x})\geq \hat{f}(\hat{x}^*)+c\|\hat{x}-\hat{x}^*\|^2.
\end{equation}
Shrinking $ \widehat{\mathcal{ N}}$ until it is contained in $\varphi(\mathcal{U})$. Define $d_1(\cdot):=d(\cdot,x^*)$, we know that $d_1(\cdot)$ is a smooth function on a neighborhood $\mathcal{U}_1$ of $x^*$ in $\mathcal{ M }$. Combining with the fact that $\varphi(\cdot)$ is a homomorphism in a neighborhood $\mathcal{U}_2\subseteq \mathcal{U}_1$ of $x^*$, there exists $\alpha>0$, such that
\begin{equation}\label{eq:man-varphi-lip}
	d_1(x) = d_1\circ \varphi^{-1}(\varphi(x))\leq \alpha\|\varphi(x)-\varphi(x^*)\|\quad\forall\,x\in \mathcal{U}_2.
\end{equation}
Combining (\ref{eq:rn-2grow}) and (\ref{eq:man-varphi-lip}), we obtain that for any $x\in \varphi^{-1}(\widehat{\mathcal{ N}})\cap \mathcal{U}_2$,
\begin{equation}\label{eq:man-2grow}
	f(x)\geq f(x^*)+\dfrac{c}{\alpha^2} d^2(x,x^*).
\end{equation}
It is clear that $\varphi^{-1}(\widehat{\mathcal{ N}})\cap \mathcal{U}_2$ is a neighborhood of $x^*$. Hence, we complete the proof.\qed

\begin{remark}\label{rm:proof-sosc}
Another approach to proving Theorem \ref{thm:man-sosc} involves an intrinsic method that does not rely on a coordinate chart. This approach is analogous to the proof provided in \cite[Theorem 3.86]{BS13}. This intrinsic method shows that the neighborhood of $ x^* $ where the quadratic growth condition holds is chart-independent. However, we use a chart-dependent approach to establish the equivalence between the M-SOSC for the manifold problem (\ref{eq:man-cone}) and the standard SOSC for the Euclidean problem (\ref{eq:rn-cone}).
\end{remark} 

\begin{remark}\label{rm:equi-sosc-nlp}
	It is worth noting that if the set $\mathcal{K}$ in (\ref{eq:prime}) is the non-negative orthant, then the manifold second-order necessary condition \eqref{eq:man-sonc} and the manifold second-order sufficient condition \eqref{eq:man-sosc} {\color{black} reduces to} those proposed in \cite[Theorem 4.2]{YZS14} and \cite[Theorem 4.3]{YZS14}, respectively, since the support function $\sigma\left(\cdot, \mathcal{T}_{\mathcal{K}}^{2}\left(g(x^*), Dg(x^*) \xi\right)\right)$ in (\ref{eq:man-sonc}) and (\ref{eq:man-sosc}) is zero (see \cite[page 177]{BS13} for details).
\end{remark}

\begin{remark}\label{rm:em-sosc-equi}
	If $\mathcal{M}$ is an embedded submanifold $\mathcal{M}=\{h(x)=0\}$ with $h'(x)$ having full row rank, then for all $\xi \in T_x\mathcal{M}$, $\left\langle \xi, \operatorname{Hess}_{x} L(x ; y) \xi\right\rangle = \left\langle \xi, \nabla^2_{xx} \tilde{L}(x ; y,z) \xi\right\rangle$, where $\tilde{L}$ and $z$ are defined in (\ref{eq:em-kkt}). Combined with the fact that the critical cones for these two problems are equal and the support function $\sigma\left(\cdot, \mathcal{T}_{\{0\}}^{2}\left(h(x), h'(x) \xi\right)\right)$ for the manifold constraint is zero since the constraint is polyhedral, the traditional Euclidean SOSC holds at $x$ if and only if the M-SOSC (\ref{eq:man-sosc}) holds at $x$.
\end{remark}

\subsection{ The equivalent perturbed problem}\label{subsec:equi-perturb}

Before characterizing the robust isolated calmness property of the KKT solution for problem~(\ref{eq:prime}), we first define the perturbed optimization problem by employing the normal coordinate chart of the Riemannian manifold ${\cal M}$. 

For a given point $x\in\mathcal{M}$, let $(\mathcal{ U},\varphi)$ be the normal coordinate chart defined around $x$. Then we are able to transform problem (\ref{eq:prime}) into the optimization problem on $\mathbb{R}^n$ around $x$  locally and define the corresponding perturbed optimization problem by:
\begin{equation}\label{eq:rn-perturb}
	\begin{array}{ll}
		\min & \hat{f}(\hat{x})+\theta(\hat{g_1}(\hat{x})+b)-\langle \hat{a}, \hat{x}\rangle \\
		\text { s.t. } & \hat{g_2}(\hat{x})+c \in \mathcal{Q},\\
		&\hat{x}\in\varphi(\mathcal{U}),
	\end{array}
\end{equation}
where $\hat{f}:=f\circ \varphi^{-1}$, $\hat{g_1}:=g_1\circ \varphi^{-1}$, $\hat{g_2}:=g_2\circ \varphi^{-1}$, $\hat{x}=\varphi(x)=(\varphi_1(x),\varphi_2(x),\cdots,\varphi_n(x))$, and $(\hat{a},b,c)\in \mathbb{R}^n\times \mathbb{Y}\times \mathbb{Z}$ is the perturbation parameter. Thus, we can transform (\ref{eq:rn-perturb}) into the following perturbed problem on $\mathcal{M}$: 
\begin{equation}\label{eq:man-perturb}
	\begin{array}{ll}
		\min & f(x)+\theta(g_1(x)+b)-\langle \hat{a},\varphi(x)\rangle\\
		\text { s.t. } 
		& g_2(x)+c \in \mathcal{Q},\\
		& x\in \mathcal{U}.
	\end{array}
\end{equation}
Note that problem (\ref{eq:man-perturb}) can be equivalently written as follow 
\begin{equation}\label{eq:man-epi}
	\begin{array}{ll}
		\min & f(x)+t-\langle \hat{a},\varphi(x)\rangle \\
		\text { s.t. } 
		& (g_1(x)+b,t)\in \operatorname{epi}\theta,\\
		& g_2(x)+c \in \mathcal{Q},\\
		& x\in \mathcal{U},
	\end{array}
\end{equation}
and problem (\ref{eq:rn-perturb}) can be correspondingly written locally as
\begin{equation}\label{eq:rn-epi}
	\begin{array}{ll}
		\min & \hat{f}(\hat{x})+t-\langle \hat{a},\hat{x}\rangle \\
		\text { s.t. } 
		& (\hat{g}_1(\hat{x})+b,t)\in \operatorname{epi}\theta,\\
		& \hat{g}_2(\hat{x})+c \in \mathcal{Q},\\
		& \hat{x}\in \varphi(\mathcal{U}).
	\end{array}
\end{equation}

Let $\widehat{\mathcal{F}}^{e}(\hat{a},b,c)$ be the set of all feasible points of the epigraph form problem (\ref{eq:rn-epi}) with a given $(\hat{a},b,c)$, 
and $\mathcal{F}^{e}(\hat{a},b,c)$ be the set of all feasible points of the epigraph form problem (\ref{eq:man-epi}) with $(\hat{a},b,c)$.
The Lagrangian function of problem (\ref{eq:man-epi}) with $(\hat{a},b,c)=(0,0,0)$ is defined by
\begin{equation}\label{eq:man-lagran-epi}
	L^{e}(x,t ; y,z,\tau):=f(x)+t+\langle y, g_1(x)\rangle+\langle z, g_2(x)\rangle+t\tau, \quad(x,t,y,z,\tau) \in \mathcal{M} \times \mathbb{R} \times \mathbb{Y} \times \mathbb{Z}  \times \mathbb{R},
\end{equation}
and the Lagrangian function of problem (\ref{eq:rn-epi}) with $(\hat{a},b,c)=(0,0,0)$ is defined by
\begin{equation}\label{eq:rn-lagran-epi}
	\widehat{L}^{e}(\hat{x},t; y,z,\tau):=\hat{f}(\hat{x})+t+\langle y, \hat{g}_1(\hat{x})\rangle+\langle z, \hat{g}_2(\hat{x})\rangle+t\tau, \quad(\hat{x},t,y,z,\tau) \in \mathbb{R}^n  \times \mathbb{R} \times \mathbb{Y} \times \mathbb{Z}  \times \mathbb{R}.
\end{equation}
For a given perturbation parameter $(\hat{a},b,c)$, the Karush–Kuhn–Tucker (KKT) optimality condition for problem (\ref{eq:rn-epi}) takes the following form:
\begin{equation}\label{eq:rn-kkt-epi}
	\left\{\begin{array} { l } 
		{ \hat{a} = \nabla _ { \hat{x} } \widehat{L}^{e} ( \hat{x},t; y,z,\tau )=\nabla \hat{f}(\hat{x})+\hat{g}'_1(\hat{x})^*y+\hat{g}'_2(\hat{x})^*z , } \\
		{0=\nabla_t \widehat{L}^{e} ( \hat{x},t; y,z,\tau )=1+\tau,}\\
		{ (y,\tau) \in \mathcal{N}_{\operatorname{epi}\theta}(\hat{g}_1 ( \hat{x} )+b,t),}\\
		{ z \in \mathcal{N}_{\mathcal{Q}}(\hat{g}_2 ( \hat{x} )+c)}.
	\end{array} \right.
\end{equation}
By (\ref{eq:co-ex}), we know that  $ \nabla _ { \hat{x} } \widehat{L}^{e} ( \hat{x},t ; y,z,\tau )=G_{\varphi(x)}D\varphi(x)\operatorname{grad}_xL^{e}(x,t ; y,z,\tau )$. Hence, from (\ref{eq:rn-kkt-epi}) we obtain the following system:
\begin{equation}\label{eq:man-kkt-epi}
	\left\{\begin{array} { l } 
		{(D\varphi(x))^{-1}G_{\varphi(x)}^{-1}\hat{a}=\operatorname{grad}_x L^{e}(x,t; y,z,\tau) =\operatorname{grad} f(x)+Dg_1(x)^*y+Dg_2(x)^*z , } \\
		{0=\nabla_t L^{e}( x,t; y,z,\tau )=1+\tau,}\\
		{ (y,\tau) \in \mathcal{N}_{\operatorname{epi}\theta}(g_1 ( x )+b,t),}\\
		{ z \in \mathcal{N}_{\mathcal{Q}}(g_2 (x )+c)}.
	\end{array}\right.
\end{equation}
This is actually the KKT system of problem (\ref{eq:man-epi}), since $(D\varphi(x))^{-1}G_{\varphi(x)}^{-1}\hat{a}=\sum_{i=1}^{n}\hat{a}_i\operatorname{grad} \varphi_i(x)$ if we write $\hat{a}=(\hat{a}_1,\hat{a}_2,\cdots,\hat{a}_n)\in\mathbb{R}^n$. Indeed, from (\ref{eq:inner-pro2}), we know that for any $\xi\in T_x\mathcal{ M}$,
$$
\begin{aligned}
	\langle (D\varphi(x))^{-1}G_{\varphi(x)}^{-1}\hat{a},\xi\rangle&= (D\varphi(x)(D\varphi(x))^{-1}G_{\varphi(x)}^{-1}\hat{a})^{T}G_{\varphi(x)}D\varphi(x)\xi=\hat{a}^{T}D\varphi(x)\xi\\
	&=\left\langle \sum_{i=1}^{n}\hat{a}_{i}\operatorname{grad} \varphi_i(x),\xi\right\rangle.
\end{aligned}
$$

For a given $(\hat{a},b,c)\in\mathbb{R}^n\times\mathbb{Y}\times\mathbb{Z},$ the set of all solutions $(\hat{x},t; y,z,\tau)$ to (\ref{eq:rn-kkt-epi})
is denoted by $\widehat{\mathcal{S}}^{e}_{\mathrm{KKT}}(\hat{a},b,c)$. 
The set of Lagrangian multipliers with $(\hat{x},t,\hat{a},b,c)$ is defined by
$$
\widehat{M}^{e}(\hat{x},t, \hat{a}, b,c):=\left\{(y,z,\tau )\in\mathbb{Y}\times\mathbb{Z}\times\mathbb{R} \mid(\hat{x},t, y,z,\tau) \in \widehat{\mathcal{S}}^{e}_{\mathrm{KKT}}(\hat{a}, b,c)\right\}.
$$
We can also define the solution set of (\ref{eq:man-kkt-epi}) as $\mathcal{S}^{e}_{\mathrm{KKT}}(\hat{a}, b,c)$ and denote the multiplier set as
$$
M^{e}(x,t, \hat{a}, b,c):=\left\{(y,z,\tau )\in \mathbb{Y}\times\mathbb{Z}\times\mathbb{R} \mid(x,t, y,z,\tau) \in \mathcal{S}^{e}_{\mathrm{KKT}}(\hat{a}, b,c)\right\}.
$$
Clearly, there are some relationships between the solution set or multiplier set of (\ref{eq:rn-kkt-epi}) and (\ref{eq:man-kkt-epi}) as below:
\begin{subequations}
	\begin{equation}\label{eq:eqv-skkt}
		\widehat{\mathcal{S}}^{e}_{\mathrm{KKT}}(\hat{a}, b,c)=\tilde{\varphi}(\mathcal{S}^{e}_{\mathrm{KKT}}(\hat{a}, b,c)\cap
		(\mathcal{U}\times\mathbb{R}\times\mathbb{Y}\times\mathbb{Z}\times\mathbb{R}))\left( \tilde{\varphi}:=(\varphi,\operatorname{id})\right) ,
	\end{equation}
	\begin{equation}\label{eq:eqv-multi}
		\widehat{M}^{e}(\hat{x},t, \hat{a}, b,c)=M^{e}(x,t, \hat{a}, b,c).
	\end{equation}
\end{subequations}

The following proposition shows that the continuity and the isolated calmness of 	$\widehat{\mathcal{S}}^{e}_{\mathrm{KKT}}$ are inherited by $\mathcal{S}^{e}_{\mathrm{KKT}}$.

\begin{proposition}\label{pro:eqv-isocalm1}
	Consider a feasible solution $(x,t)\in\mathcal{F}^{e}(\hat{a},b,c)$. If $\widehat{\mathcal{S}}^{e}_{\mathrm{KKT}}(\hat{a}, b,c)$ is lower semicontinuous or (robustly) isolated calm at $(\hat{a},b,c)$ for $(\hat{x},t,y,z,-1)$, then  $\mathcal{S}^{e}_{\mathrm{KKT}}(\hat{a}, b,c)$ is lower semicontinuous or (robustly) isolated calm at $(\hat{a},b,c)$ for $(x,t,y,z,-1)$.
\end{proposition}
\proof For any open neighborhood $\mathcal{V}$ of $(x,t,y,z,-1)$, by (\ref{eq:eqv-skkt}), we know that $\tilde{\varphi}(\mathcal{V})$ is an open neighborhood of $(\hat{x},t,y,z,-1)$. If  $\widehat{\mathcal{S}}^{e}_{\mathrm{KKT}}(\hat{a}, b,c)$ is lower semicontinuous at $(\hat{a},b,c)$ for $(x,t,y,z,-1)$, then there exists an open neighborhood $\mathcal{ N}$ of $(\hat{a},b,c)$, such that
$$
\emptyset\neq \widehat{\mathcal{S}}^{e}_{\mathrm{KKT}}(a_1,b_1,c_1)\cap \tilde{\varphi}(\mathcal{V})\quad\forall\,(a_1,b_1,c_1)\in \mathcal{ N}.
$$
Since $\varphi$ is a homomorphism, by (\ref{eq:eqv-skkt}), we have
\begin{equation}\label{eq:eqv-skkt-tilde}
	\widehat{\mathcal{S}}^{e}_{\mathrm{KKT}}(a_1,b_1,c_1)\cap \tilde{\varphi}(\mathcal{V})\subseteq\tilde{\varphi}\left( \mathcal{S}^{e}_{\mathrm{KKT}}(a_1,b_1,c_1)\right) \cap \tilde{\varphi}(\mathcal{V})=\tilde{\varphi}\left(\mathcal{S}^{e}_{\mathrm{KKT}}(a_1,b_1,c_1)\cap \mathcal{V}\right).
\end{equation}
Thus $\emptyset\neq \mathcal{S}^{e}_{\mathrm{KKT}}(a_1,b_1,c_1)\cap \mathcal{V} ,\;\forall\,(a_1,b_1,c_1)\in \mathcal{ N}$. This means that $\mathcal{S}^{e}_{\mathrm{KKT}}(\hat{a}, b,c)$ is lower semicontinuous at $(\hat{a},b,c)$ for $(x,t,y,z,-1)$.

If $\widehat{\mathcal{S}}^{e}_{\mathrm{KKT}}(\hat{a}, b,c)$ is isolated calm at $(\hat{a},b,c)$ for $(\hat{x},t,y,z,-1)$, then there exist open neighborhoods  $\mathcal{ N}$ of $(\hat{a},b,c)$ and $\widehat{\mathcal{V}}$ of $(\hat{x},t,y,z,-1)$ and a constant $\kappa>0$, such that for any $ (\hat{x}_1,t_1,y_1,z_1,\tau_1)\in \widehat{\mathcal{S}}^{e}_{\mathrm{KKT}}(a_1,b_1,c_1)\cap \widehat{\mathcal{V}}$ and $(a_1,b_1,c_1)\in \mathcal{ N}$,
$$
\|(\hat{x}_1,t_1,y_1,z_1,\tau_1)-(\hat{x},t,y,z,-1)\|\leq \kappa\|(a_1,b_1,c_1)-(\hat{a},b,c)\|.
$$
By (\ref{eq:man-varphi-lip}), there exists $\alpha>0,$ such that for any $\tilde{x}\in\mathcal{U}$,
$d(\tilde{x},x)\leq \alpha\|\varphi(\tilde{x})-\varphi(x)\|.$
Therefore, for any $ (x_1,t_1,y_1,z_1,\tau_1)\in \mathcal{S}^{e}_{\mathrm{KKT}}(a_1,b_1,c_1)\cap
(\mathcal{U}\times\mathbb{R}\times\mathbb{Y}\times\mathbb{Z}\times\mathbb{R})\cap \tilde{\varphi}^{-1}(\widehat{\mathcal{V}})$ and $(a_1,b_1,c_1)\in \mathcal{ N}$,
$$
d\left( (x_1,t_1,y_1,z_1,\tau_1),(x,t,y,z,-1)\right) \leq (\alpha+1)\kappa\|(a_1,b_1,c_1)-(\hat{a},b,c)\|.
$$
This implies that  $\mathcal{S}^{e}_{\mathrm{KKT}}(\hat{a}, b,c)$ is isolated calm at $(\hat{a},b,c)$ for $(x,t,y,z,-1)$. 

Moreover, if for each $(a_1,b_1,c_1)\in \mathcal{ N}$,  $\tilde{\varphi}\left(\mathcal{S}^{e}_{\mathrm{KKT}}(a_1,b_1,c_1)\cap
(\mathcal{U}\times\mathbb{R}\times\mathbb{Y}\times\mathbb{Z}\times\mathbb{R})\cap \tilde{\varphi}^{-1}(\widehat{\mathcal{V}})\right)\neq\emptyset$, then $\mathcal{S}^{e}_{\mathrm{KKT}}(a_1,b_1,c_1)\cap \mathcal{V}\neq\emptyset$. Thus, $\mathcal{S}^{e}_{\mathrm{KKT}}(\hat{a}, b,c)$ is robustly isolated calm at $(\hat{a},b,c)$ for $(x,t,y,z,-1).$\qed

\subsection{The characterization of the robust isolated calmness of the KKT solution mapping}\label{subsec:isolated}
In this subsection, we will characterize the robust isolated calmness property of the KKT solution mapping for problem (\ref{eq:man-perturb}). We first present the analysis of the perturbation properties of $\mathcal{S}^{e}_{\mathrm{KKT}}(\hat{a}, b,c)$  for (\ref{eq:man-epi}). Let $(\hat{a},b,c)=(0,0,0)$ in (\ref{eq:man-epi}), recalling Definition \ref{def:rcq+srcq}, we say the M-RCQ  holds at a feasible solution $(x,\theta(g_1(x)))$ if 
\begin{equation}\label{eq:man-rcq-epi}
	\left[\begin{array}{c}
		\left(Dg_1(x), 1\right)\\
		\left(Dg_2(x), 0\right)
	\end{array}\right](T_{x}\mathcal{ M} \times \mathbb{R})+\left[\begin{array}{c}
		\mathcal{T}_{\operatorname{epi}\theta}(g_1(x),\theta(g_1(x)))\\
		\mathcal{T}_{\mathcal{Q}}(g_2(x))
	\end{array}\right]=\left[\begin{array}{c}
		\mathbb{Y} \times \mathbb{R}\\
		\mathbb{Z} \\
	\end{array}\right],
\end{equation}
and say the M-SRCQ  holds at $(x,\theta(g_1(x)))$ with respect to $(y, z,-1) \in M^{e}(x, \theta(g_1(x)), 0,0,0)$ if
\begin{equation}\label{eq:man-srcq-epi}
	\left[\begin{array}{c}
		\left(Dg_1(x), 1\right)\\
		\left(Dg_2(x), 0\right)
	\end{array}\right](T_{x}\mathcal{ M} \times \mathbb{R})+\left[\begin{array}{c}
		\mathcal{T}_{\operatorname{epi}\theta}(g_1(x),\theta(g_1(x)))\cap(y,-1)^{\perp}\\
		\mathcal{T}_{\mathcal{Q}}(g_2(x))\cap z^{\perp}
	\end{array}\right]=\left[\begin{array}{c}
		\mathbb{Y} \times \mathbb{R}\\
		\mathbb{Z} \\
	\end{array}\right].
\end{equation}
The critical cone at a feasible point $(x,\theta(g_1(x)))$ for problem (\ref{eq:man-epi}) takes the form of
\begin{equation}\label{eq:man-cri-epi}
	\begin{array}{ll}
		\mathcal{C}^{e}(x, \theta(g_1(x)))
		:=\left\{\left(d_{1}, d_{2}\right) \in T_{x}\mathcal{M} \times \mathbb{R}\right. \mid & \left(Dg_1(x)d_{1}, d_{2}\right) \in \mathcal{T}_{\operatorname{epi}\theta}\left( g_1(x), \theta(g_1(x))\right) ,   \\
		&\left.  Dg_2(x) d_{1} \in \mathcal{T}_{\mathcal{Q}}(g_2(x)),\  Df(x) d_{1}+d_{2} \leq 0\right\}.
	\end{array}
\end{equation}
If $x$ is further a stationary point for (\ref{eq:rn-epi}) with $(\hat{a},b,c)=(0,0,0)$ and there exists $(y,z,-1)\in M^{e}(x,\theta(g_1(x)),0,0,0)$, then the critical cone can be written as
\begin{equation}\label{eq:man-cri2-epi}
	\begin{array}{ll}
		&\mathcal{C}^{e}(x, \theta(g_1(x)))\\
		= &\left\{\left(d_{1}, d_{2}\right) \mid 
		\left(Dg_1(x)d_{1}, d_{2}\right) \in \mathcal{T}_{\operatorname{epi}\theta}\left( g_1(x), \theta(g_1(x))\right), \ Dg_2(x) d_{1} \in \mathcal{T}_{\mathcal{Q}}(g_2(x)), \  Df(x) d_{1}+d_{2}= 0\right\}\\
		= &\left\{\left(d_{1}, d_{2}\right) \mid 
		\left(Dg_1(x)d_{1}, d_{2}\right) \in \mathcal{T}_{\operatorname{epi}\theta}\left( g_1(x), \theta(g_1(x))\right)\cap(y,-1)^{\perp},\ Dg_2(x) d_{1} \in \mathcal{T}_{\mathcal{Q}}(g_2(x))\cap z^{\perp}\right\}.
	\end{array}
\end{equation}
From now on, we always assume that the function $\theta$ is $\mathcal{C}^2$-cone reducible at $g_1(x)$ and the set $\mathcal{Q}$ is $\mathcal{C}^2$-cone reducible at $g_2(x)$ (Definition \ref{def:c2reducible}). Then, we know that the M-SOSC at  $(x,\theta(g_1(x)))$ for problem (\ref{eq:man-epi}) with $(\hat{a},b,c)=(0,0,0)$ holds if for any $d \in \mathcal{C}^{e}(x, \theta(g_1(x))) \backslash\{0\}$,
\begin{equation}\label{eq:man-sosc-epi}
	\begin{aligned}
		\sup _{(y, z,-1) \in M^{e}}&\left\langle d, \operatorname{Hess}_{(x, t)} L^{e}(x, \theta(g_1(x)), y, z,-1) d\right\rangle\\
		&-\sigma\left((y,-1,z), \mathcal{T}_{\operatorname{epi}\theta\times \mathcal{Q}}^{2}\left( \left( g_1(x)  , \theta(g_1(x)),g_2(x)\right)  ,\left(\left( Dg_1(x), 1\right)d ,Dg_2(x) d\right) \right)\right) >0.
	\end{aligned}
\end{equation}

Now let us return to the original perturbed problem (\ref{eq:man-perturb}). Let $(\hat{a},b,c)$ be given. We say that $x$ is a feasible solution to problem (\ref{eq:man-perturb}) if
\begin{equation}\label{eq:man-feas-perturb}
	x\in \mathcal{F}(\hat{a},b,c):=\{g_1(x)+b\in \operatorname{dom}\theta\mid g_2(x)+c\in\mathcal{Q}\}.
\end{equation}
Denote $l:\mathcal{ M}\times \mathbb{Z}\rightarrow \mathbb{R}$ by
\begin{equation}\label{eq:man-lagrang-peturb}
	l(x,z):=f(x)+\langle z,g_2(x)\rangle,\quad(x,z)\in \mathcal{ M}\times \mathbb{Z},
\end{equation}
then the KKT condition takes the form of 
\begin{equation}\label{eq:man-kkt-perturb1}
	\left\{\begin{array} { l } 
		{(D\varphi(x))^{-1}G_{\varphi(x)}^{-1}\hat{a}\in\operatorname{grad} _xl(x;z) +\partial\theta\circ( g_1(x)+b) , } \\
		{ z \in \mathcal{N}_{\mathcal{Q}}(g_2 (x )+c)}.
	\end{array}\right.
\end{equation}
By Proposition \ref{pro:subdiff} and \cite[Proposition 2.3.6]{C90}, we know that $\partial \theta\circ \left( g_1(x)+b\right)=Dg_1(x)^*\partial\theta (g_1(x)+b)$. Therefore, the KKT condition (\ref{eq:man-kkt-perturb1}) can be rewritten as 
\begin{equation}\label{eq:man-kkt-perturb2}
	\left\{\begin{array} { l } 
		{(D\varphi(x))^{-1}G_{\varphi(x)}^{-1}\hat{a}=\operatorname{grad}  f(x)+Dg_1(x)^*y+Dg_2(x)^*z , } \\
		{y\in \partial \theta(g_1(x)+b),}\\
		{ z \in \mathcal{N}_{\mathcal{Q}}(g_2 (x )+c)}.
	\end{array}\right.
\end{equation}	
Let $\mathcal{S}_{\mathrm{KKT}}:\mathbb{R}^n\times\mathbb{Y}\times\mathbb{Z}\rightarrow \mathcal{ M}\times\mathbb{Y}\times \mathbb{Z}$ be the following KKT solution mapping: 
\begin{equation}\label{eq:man-skkt-perturb}
	\mathcal{S}_{\mathrm{KKT}}(\hat{a},b,c)=\{(x,y,z)\in \mathcal{ M}\times \mathbb{Y}\times \mathbb{Z}\mid (x,y,z)\text{ satisfies }(\ref{eq:man-kkt-perturb2})\}.
\end{equation}
We also denote $M(x,\hat{a},b,c)$ as the set of Lagrange multipliers with respect to $x$.

The following proposition establishes the relation between $\mathcal{S}^{e}_{\mathrm{KKT}}$ and $\mathcal{S}_{\mathrm{KKT}}$. The proof is similar to that of Proposition 3.2 in \cite{CS17}, and we include the proof here for completeness.
\begin{proposition}\label{pro:eqv-isocalm2}
Let $(x^*, \theta(g_1(x^*))) \in \mathcal{M} \times \mathbb{R}$ be a local optimal solution to problem (\ref{eq:man-epi}), with $M^{e}(x^*, \theta(g_1(x^*)), 0, 0, 0) \neq \emptyset$. Assume $(y^*, z^*, -1) \in M^{e}(x^*, \theta(g_1(x^*)), 0, 0, 0)$. Then the KKT solution mapping $\mathcal{S}^{e}_{\mathrm{KKT}}$ is robustly isolated calm at the origin for $(x^*, \theta(g_1(x^*)), y^*, z^*, -1)$ if and only if the KKT solution mapping $\mathcal{S}_{\mathrm{KKT}}$ is robustly isolated calm at the origin for $(x^*, y^*, z^*)$.
\end{proposition}
\proof  Given any $(\hat{a}, b, c) \in \mathbb{R}^n \times \mathbb{Y} \times \mathbb{Z}$ and $(x, y, z) \in \mathcal{S}_{\mathrm{KKT}}(\hat{a}, b, c)$, we are able to show that $(x, \theta(g_1(x)), y, z, -1) \in \mathcal{S}^e_{\mathrm{KKT}}(\hat{a}, b, c)$. Since $g_1 \in \mathcal{C}^2$ and $\theta$ is convex, there exists a constant $\kappa_1 > 0$ such that for any $x$ sufficiently close to $x^*$,
\[
\left\| \theta(g_1(x)) - \theta(g_1(x^*)) \right\| \leq \kappa_1 d(x, x^*).
\]
This implies that the robustly isolated calmness of $\mathcal{S}_{\mathrm{KKT}}$ at $(x^*, y^*, z^*)$ ensures the robustly isolated calmness of $\mathcal{S}^e_{\mathrm{KKT}}$ at $(x^*, \theta(g_1(x^*)), y^*, z^*, -1)$.

Conversely, using \cite[Corollary 2.4.9]{C90}, we have:
\[
(y, -1) \in \mathcal{N}_{\operatorname{epi} \theta}(g_1(x) + b, \theta(g_1(x) + b)) \Leftrightarrow y \in \partial \theta(g_1(x) + b).
\]
Thus, if $(x, \theta(g_1(x)), y, z, -1) \in \mathcal{S}^e_{\mathrm{KKT}}(\hat{a}, b, c)$, it follows that $(x, y, z) \in \mathcal{S}_{\mathrm{KKT}}(\hat{a}, b, c)$. The ``only if" part directly follows from the definition of robust isolated calmness.
\qed

By noting that $\theta$ is convex on $\mathbb{Y}$, we know from \cite[Proposition 8.12, Exercise 8.4]{RW09} and \cite[Proposition 2.2.7]{C90} that 
$$
\partial \theta(g_1(x)+b)= \{y\mid\langle y,d\rangle\leq\theta^{\downarrow}\left( g(x),d\right),\text{ for all }d\in\mathbb{Y} \}.
$$
Thus, the KKT condition (\ref{eq:man-kkt-perturb2}) can be equivalently written as 
\begin{equation}\label{eq:kkt-equi-direc}
	\left\{\begin{array} { l } 
		\operatorname{grad} _xl(x;z)-(D\varphi(x))^{-1}G_{\varphi(x)}^{-1}\hat{a} =-Dg_1(x)^{*}y,\\
		\theta^{\downarrow}\left( g_1(x)+b,d\right) -\langle y,d\rangle\geq 0  \quad \forall\,d\in \mathbb{Y},\\
		z \in \mathcal{N}_{\mathcal{Q}}(g_2 (x )+c).
	\end{array}\right.
\end{equation}
Let $(\hat{a},b,c)=(0,0,0)$. Define the critical cone of the function $\theta$ and $g$ by 
\begin{equation}\label{eq:critical-theta}
	\mathcal{C}_{\theta,g}(x,y):=\left\{d \in \mathbb{Y} \mid \theta^{\downarrow}(g(x) ; d)=\langle d, y\rangle \right\}.
\end{equation}
The M-RCQ is said to hold at a feasible solution $x$ of the problem (\ref{eq:man-perturb}) if
\begin{equation}\label{eq:man-rcq-perturb}
	\left[\begin{array}{c}
		Dg_1(x)\\
		Dg_2(x)
	\end{array}\right]T_{x}\mathcal{ M}+\left[\begin{array}{c}
		\mathcal{T}_{\operatorname{dom}\theta}(g_1(x))\\
		\mathcal{T}_{\mathcal{Q}}(g_2(x))
	\end{array}\right]=\left[\begin{array}{c}
		\mathbb{Y} \\
		\mathbb{Z} \\
	\end{array}\right],
\end{equation}
and by (\ref{eq:tangent-epideri}), the M-SRCQ is said to hold at a stationary point $x$ with respect to  $(y,z)\in M(x,0,0,0)$ if
\begin{equation}\label{eq:man-srcq-perturb}
	\left[\begin{array}{c}
		Dg_1(x)\\
		Dg_2(x)
	\end{array}\right]T_{x}\mathcal{ M}+\left[\begin{array}{c}
		\mathcal{C}_{\theta,g_1}(x,y)\\
		\mathcal{T}_{\mathcal{Q}}(g_2(x))\cap z^\perp
	\end{array}\right]=\left[\begin{array}{c}
		\mathbb{Y} \\
		\mathbb{Z} \\
	\end{array}\right].
\end{equation}
The critical cone at a stationary point $x$ of problem (\ref{eq:man-perturb}) with $z\in M(x,0,0,0)$ is
\begin{equation}\label{eq:man-cri-perturb}
	\mathcal{C}(x)=\left\{\xi \in T_{x}\mathcal{M} \mid Dg_2(x)\xi \in \mathcal{T}_{\mathcal{Q}}(g_2(x))\cap z^{\perp}, Df(x)\xi+(\theta\circ g_1)^{\circ}(x ; \xi) = 0\right\}.
\end{equation}
By (\ref{eq:2d-tan}) and the $\mathcal{C}^{2}$-cone reducibility of $\theta$ and $\mathcal{ Q}$, we know that  the M-SOSC at $x$  is said to hold if for any $\xi \in \mathcal{C}(x) \backslash\{0\}$,
\begin{equation}\label{eq:man-sosc-perturb}
	\sup _{y,z \in M(x,0, 0,0)}\left\lbrace \left\langle \xi, \operatorname{Hess}_{x} l(x, z)\xi\right\rangle-\psi_{(g_1(x),Dg_1(x)\xi)}^{*}(y)-\sigma\left(z, \mathcal{T}_{\mathcal{Q}}^{2}\left( g_2(x) ; Dg_2(x) \xi\right) \right)\right\rbrace >0,
\end{equation}
where $\psi_{(g_1(x),Dg_1(x)\xi)}^{*}(\cdot)$ is the conjugate function of $\psi_{(g_1(x),Dg_1(x)\xi)}(\cdot)=\theta_{-}^{\downarrow\downarrow}(g_1(x) ; Dg_1(x)\xi, \cdot)$ for any $g_1(x) \in \operatorname{dom} \theta$ and any $\xi \in T_x\mathcal{M}$.

\begin{remark}\label{rm:example-sosc-srcq}
Now let us consider the M-SRCQ and M-SOSC conditions for problem (\ref{exam:sphere}), where $f(x)=x_2^2$, $g_1(x)=x_1-x_2$, $\theta(\cdot)=|\cdot|$ and $g_2(x)=2x_1+x_2$. We focus on the properties at the global solution $x^{*}=\left({\sqrt{2} }/{2},{\sqrt{2}}/{2}\right)^{\top}$ with respect to the multipliers $y^*={\sqrt{2}}/{2}$ and $z^*=0$. The tangent space of $\mathcal{S}^{1}=\{x=(x_1,x_2)\mid x_1^2+x_2^2=1\}$ at $x^*$ is $T_{x^*}\mathcal{S}^{1}=\left\lbrace \xi\in \mathbb{R}^{2}\mid x^{*\top}\xi=0 \right\rbrace =\{ \left( v,-v \right)^{\top} \mid v\in \mathbb{R}  \}$. 
	Thus, for any $ u\in Dg_1(x^*)T_{x^*}\mathcal{S}^1,$ there exists $(v,-v)^{\top}\in \mathbb{R}^2,$ such that 
	$$u= \langle \nabla g_1(x^*),(v,-v)^{\top}\rangle = \big\langle (1,-1)^{\top}, (v,-v)^{\top}\big\rangle =2v. $$
	Since $v$ can be taken arbitrary, $ Dg_1(x^*)T_{x^*}\mathcal{S}^1=\mathbb{R}$. Combining with the fact that $\mathcal{T}_{\mathcal{Q}}\left( g_2(x^*)\right) =\mathbb{R }$ and $(z^{*})^{\perp}=(0)^{\perp}=\mathbb{R}$, the M-SRCQ condition is satisfied at $x^*$ with respect to $(y^*,z^*)$.
	Moreover, using \cite[Lemma 5.1 and Definition 5.1]{YZS14} and the fact that $\theta\circ g_1(\cdot)$ is convex on $\mathbb{R }^2$, for any $\xi=\left( v,-v\right) \in T_{x^*}S^{1}$, 
	$$
	\begin{aligned}
	\left( \theta\circ g_1\right)^{\circ}(x^*,\xi)&=\left( \theta\circ g_1\right)'(x^*,\xi)=\lim_{t\downarrow 0}\dfrac{\left| x^*_1+tv-(x^*_2-tv)\right| -\left| x^*_1-x^*_2\right| }{t}=2\left| v\right|.
	\end{aligned}
	$$
	Therefore, the critical cone at $x^*$ is
	$$
	\begin{aligned}
	\mathcal{C}(x^*)&=\left\lbrace \xi=\left\langle v,-v\right\rangle \mid Dg_2(x^*)\xi\in \mathbb{R}, Df(x^*)\xi+ \left( \theta\circ g_1\right) ^{\circ}(x^*,\xi)= 0\right\rbrace \\
	&=\left\lbrace \left\langle v,-v\right\rangle\mid -\sqrt{2}v+2\left| v \right|= 0 \right\rbrace =\left\lbrace (0,0)\right\rbrace ,
	\end{aligned}
	$$
	and the M-SOSC condition is thus satisfied.
\end{remark}

Next, we shall study the characterization of the robust isolated calmness of $\mathcal{S}_{\mathrm{KKT}}$ at $(0,0,0)$ for a local optimal solution $x^*$ of problem (\ref{eq:man-perturb}). First, let $(\hat{x},t) \in \mathbb{R}^n$ be a feasible solution to problem (\ref{eq:rn-epi}) with $(\hat{a}, b,c)=(0,0,0)$.  It is known (cf. \cite[Theorem 24]{DSZ17} and \cite[Proposition 3.1]{CS17}) that the  mapping $\widehat{\mathcal{ S}}^{e}_\mathrm{KKT}$ has the isolated calmness property at $(0,0,0)$ for $(\hat{x},t)$ if the SOSC holds at $(\hat{x},t)$ and SRCQ holds at $(\hat{x},t)$ with respect to the multiplier. Therefore, the following results for the semicontinuity and isolated calmness of $\mathcal{ S}^{e}_\mathrm{KKT}$ of the problem (\ref{eq:man-epi}) extends \cite{DSZ17} to manifold optimization, which can be easily obtained by using the normal coordinate chart around $x^*$. Combining Proposition \ref{pro:eqv-isocalm1} with \cite[Theorem 17]{DSZ17}, we obtain the following theorem. 

\begin{theorem}\label{thm:sosc+srcq-epi-lsc}
	Let $(x^*,\theta(g_1(x^*)))$ be a feasible solution of problem (\ref{eq:man-epi}) with $(\hat{a}, b,c) = (0, 0,0)$. Suppose that the M-SRCQ (\ref{eq:man-srcq-epi}) holds at $(x^*,\theta(g_1(x^*)))$ with respect to $(y^*, z^*,-1) \in M^{e}(x^*,$ $ \theta(g_1(x^*)), 0,0,0)$ and the M-SOSC (\ref{eq:man-sosc-epi}) holds at $(x^*,\theta(g_1(x^*)))$ for problem (\ref{eq:man-epi}) with respect to $(\hat{a},b,c) = (0,0,0)$. Then the set-valued mapping $\mathcal{S}^{e}_{\mathrm{KKT}}$ is lower semicontinuous at $(0,0, 0, x^*,\theta(g_1(x^*)),y^*,z^*,-1) \in \operatorname{gph} \mathcal{S}^{e}_{\mathrm{KKT}}$.
\end{theorem}

By applying \cite[Theorem 24]{DSZ17} to $\widehat{S}^{e}_{\mathrm{KKT}}$ and using Proposition \ref{pro:eqv-isocalm1}, we obtain the following result on the characterization of the (robust) isolated calmness of $\mathcal{S}^{e}_{\mathrm{KKT}}$.
\begin{theorem}\label{thm:man-srcq+sosc-epi}
	Let $(x^*,\theta(g_1(x^*)))$ be a feasible solution to problem (\ref{eq:man-epi}) with $(\hat{a}, b,c)=(0,0,0)$. Suppose that the M-RCQ (\ref{eq:man-rcq-epi}) holds at $(x^*,\theta(g_1(x^*))).$ Let $(y^*, z^*,-1) \in M^{e}(x^*, \theta(g_1(x^*)), 0,0,0) \neq \emptyset$. Then the following statements are equivalent:
	\begin{itemize}
		\item[(i)] the M-SRCQ (\ref{eq:man-srcq-epi}) holds at $(x^*,\theta(g_1(x^*)))$ with respect to $(y^*, z^*,-1)$ and the M-SOSC (\ref{eq:man-sosc-epi}) holds at $(x^*,\theta(g_1(x^*)))$ for problem (\ref{eq:man-epi}) with $(\hat{a}, b,c) = (0, 0,0)$;
		\item[(ii)] $(x^*,\theta(g_1(x^*)))$ is a locally optimal solution to problem (\ref{eq:man-epi}) with $(\hat{a}, b,c)=(0,0,0)$ and $\mathcal{S}^{e}_{\mathrm{KKT}}$ is (robustly) isolated calm at the origin for $(x^*,\theta(g_1(x^*)), y^*, z^*,-1)$.
	\end{itemize}
\end{theorem}

When $(\hat{a}, b,c) = (0,0, 0)$, the KKT system (\ref{eq:man-kkt-perturb2}) is equivalent to the following system of nonsmooth equations:
$	0\in F(x, y,z),$
where $F: \mathcal{M}\times\mathbb{Y}\times\mathbb{Z}\rightarrow T\mathcal{M }\times\mathbb{Y}\times\mathbb{Z}$ is the natural mapping defined by
\begin{equation}\label{eq:man-natmap}
	F(x, y,z):=\left[\begin{array}{c}
		\operatorname{grad}  f(x)+Dg_1(x)^{*}y+Dg_2(x)^{*}z \\
		g_1(x)-\operatorname{prox}_{\theta}(g_1(x)+y)\\
		g_2(x)-\Pi_{\mathcal{Q}}(g_2(x)+z)
	\end{array}\right], \quad(x,y, z) \in \mathcal{M} \times \mathbb{Y}\times \mathbb{Z}.
\end{equation}
where $\operatorname{prox}_{\theta}(\cdot)$ is the proximal mapping of  $\theta$ which is defined as $\operatorname{prox}_{\theta}(u):=\argmin_{y\in\mathbb{Y}}\;\theta(y)+\dfrac{1}{2}\|u-y\|^2$, and $\Pi_{\mathcal{Q}}(\cdot)$ denotes the projection onto $\mathcal{Q}$. It is clear that $(0, 0,0, x^*,y^*,z^*)\in \operatorname{gph} \mathcal{S}_{\mathrm{KKT}}$ if and only if $(0,0, 0, x^*,y^*,z^*) \in \operatorname{gph} F^{-1}$. The proof of Lemma \ref{lem:man-eqv-natskkt} is almost the same as the proof of  \cite[Lemma 18]{DSZ17}, and we omit it here.

\begin{lemma}\label{lem:man-eqv-natskkt}
	Let $(0, 0,0, x^*,y^*,z^*) \in \operatorname{gph} \mathcal{S}_{\mathrm{KKT}} .$ The set-valued mapping $\mathcal{S}_{\mathrm{KKT}}$ is isolated calm at the origin for $( x^*,y^*,z^*)$ if and only if the set-valued mapping $F^{-1}$ is isolated calm at the origin for $(x^*,y^*,z^*)$.
\end{lemma}

Combining Lemma \ref{lem:man-eqv-natskkt}, Proposition \ref{pro:eqv-isocalm2} and Theorem \ref{thm:man-srcq+sosc-epi}, we can derive the main result of this section in the next Proposition.

\begin{proposition}\label{pro:man-sosc+srcq-perturb}
	Let $x^*\in\mathcal{ M }$ be a local optimal solution of problem (\ref{eq:man-perturb}) with $(\hat{a},b,c)=(0,0,0).$ Suppose that the M-RCQ (\ref{eq:man-rcq-perturb}) holds at $x^*$.  Let $(y^*,z^*)\in M(x^*,0,0,0)$. Then the following statements are equivalent: 
	\begin{itemize}
		\item [(i)]The M-SOSC (\ref{eq:man-sosc-perturb}) holds at $x^*$ and the M-SRCQ (\ref{eq:man-srcq-perturb}) holds at $x^*$ with respect to $(y^*,z^*)$;
		\item [(ii)]The KKT solution mapping $\mathcal{S}_{\mathrm{KKT}}$ is robustly isolated calm at the origin for $(x^*,y^*,z^*)$;
  	\item [(iii)]The set-valued mapping $F^{-1}$ is robustly isolated calm at the origin for $(x^*,y^*,z^*)$.
	\end{itemize}
\end{proposition}

\begin{remark}\label{rm:iso-equi-naturalmapping}
	The isolated calmness of $F^{-1}$ at the origin for $(x^*,y^*,z^*)$ implies that there exist $\kappa>0$, and a neighborhood $\mathcal{V}$ of $(x^*,y^*,z^*)$ such that
	\begin{equation}\label{eq:man-resi-bound}
		d(x,x^*)+\left\|y-y^* \right\| +\left\|z-z^* \right\| \leq \kappa \left\|F(x,y,z)\right\|\quad \forall\,(x,y,z)\in\mathcal{ V}.
	\end{equation}
\end{remark}

\section{Riemannian augmented Lagrangian method}\label{sec:alm}

In this section, we shall present the inexact augmented Lagrangian method for solving the problem (\ref{eq:prime}) and establish its local convergence analysis. When $\mathcal{ Q}$ is polyhedral, Algorithm \ref{alg:alm} will reduce to the Riemannian augmented Lagrangian method proposed in \cite{ZBDZ21}. 

The Lagrangian function of (\ref{eq:prime}) is defined by
\begin{equation}\label{eq:prim-lagrang}
	L(x, y,z):=f(x)+\langle y, g_1(x)\rangle+\langle z, g_2(x)\rangle, \quad(x,y,z) \in \mathcal{M} \times \mathbb{Y}  \times \mathbb{Z} ,
\end{equation}
and the augmented Lagrangian function $L_\rho : \mathcal{M}  \times \mathbb{Y}  \times \mathbb{Z}\rightarrow \mathbb{R}$ is in the form of 
\begin{equation}\label{eq:prim-Al}
	L_{\rho}(x, y,z):=f(x)+\theta^{\rho}\left(g_1(x)+\dfrac{y}{\rho}\right)+\frac{\rho}{2} \operatorname{dist}^{2}\left(g_2(x)+\frac{z}{\rho},\mathcal{Q}\right),
\end{equation}
where $\theta^{\rho}$ is the Moreau-Yosida regularization of $\theta$ defined by $\theta^{\rho}\left(u\right):=\min_{y\in\mathbb{Y}}\; \theta(y)+\dfrac{\rho}{2}\|u-y\|^2$. 
Moreover, we define the auxiliary function
\begin{equation}\label{eq:auxiliary}
	V(x, y,z, \rho):=\max\left\{\left\|g_1(x)-\operatorname{prox}_{\theta}\left( g_1(x)+\dfrac{y}{\rho}\right) \right\|,\left\|g_2(x)-\Pi_{\mathcal{Q}}\left(g_2(x)+\dfrac{z}{\rho}\right)\right\|\right\} .
\end{equation}
 The detail of Riemannian ALM for solving (\ref{eq:prime}) is presented in Algorithm \ref{alg:alm}.

\begin{algorithm}[h]
	\caption{Riemannian augmented Lagrangian method}  
	\label{alg:alm}
	\hspace*{0.02in} {\bf Input:}  
	Let $\left(x^{0},y^{0},z^{0}\right) \in \mathcal{M} \times \mathbb{Y}\times \mathbb{Z}$, $B \subseteq \mathbb{Y}\times\mathbb{Z}$ be bounded, $ \rho^{0}>0$, $\gamma>1$, $\alpha \in(0,1)$, a sequence $\{\epsilon_{k}\}\in\mathbb{R}_{+}$ convergence to 0, and set $k:=0$.
	\begin{algorithmic}[1]
		\State If $\left(x^{k}, y^{k},z^{k}\right)$ satisfies a suitable termination criterion: STOP.
		\State Choose $(w^{k},p^k) \in B$ and compute an $x^{k+1}$, such that
		\begin{equation}\label{eq:subproblem}
			x^{k+1}\approx \argmin_{x\in\mathcal{ M }} L_{\rho^{k}}(x,w^{k},p^{k}):=  f(x)+\theta^{\rho^k}\left(g_1(x)+\dfrac{w^k}{\rho^k}\right)+\dfrac{\rho^k}{2} \operatorname{dist}^{2}\left(g_2(x)+\dfrac{p^k}{\rho^k},\mathcal{Q}\right).
		\end{equation}
		Specially, we need to find $x^{k+1}$ satisfying
		\begin{equation}\label{epsi-grad}
			\|\operatorname{grad}  L_{\rho^{k}}(x,w^{k},p^{k})\|\leq \epsilon_k .
		\end{equation}
		\State Update the vector of multipliers to 
		\begin{equation}\label{eq:update-y}
			y^{k+1}=\rho^{k}\left[g_1(x^{k+1})+\dfrac{w^k}{\rho^{k}}-\operatorname{prox}_{\theta/\rho^{k}}\left(g_1(x^{k+1})+\dfrac{w^k}{\rho^{k}}\right)\right],
		\end{equation}
		\begin{equation}\label{eq:update-z}
			z^{k+1}=\rho^{k}\left[g_2(x^{k+1})+\dfrac{p^k}{\rho^{k}}-\Pi_{\mathcal{Q}}\left(g_2(x^{k+1})+\dfrac{p^k}{\rho^{k}}\right)\right].
		\end{equation}
		\State If $k=0$ or
		\begin{equation}\label{eq:update-rho}
			V\left(x^{k+1}, w^{k},p^k, \rho^{k}\right) \leq \tau V\left(x^{k}, w^{k-1},p^{k-1}, \rho^{k-1}\right)
		\end{equation}
		holds, set $\rho^{k+1}=\rho^{k} ;$ otherwise, set $\rho^{k+1}=\gamma \rho^{k} .$
		\State 	Set $k=k+1$ and go to step 1.
	\end{algorithmic}
\end{algorithm}

\subsection{Local convergence analysis of the Riemannian ALM}\label{subsec:alm-conver}
In this subsection, motivated by the proof in \cite{KS19}, we present the local convergence result of Algorithm~\ref{alg:alm}. We first define the KKT residual mapping as follows:
\begin{equation}\label{eq:resi-map}
	R(x,y,z):=\|\operatorname{grad}_x{L}(x;y,z)\|+\left\|g_1(x)-\operatorname{prox}_{\theta}(g_1(x)+y) \right\|
	+\left\|g_2(x)-\Pi_{\mathcal{Q}}(g_2(x)+z)\right\|.
\end{equation}

The following theorem provides a characterization of the local error bound for the distance between points near a KKT solution $(x^*, y^*, z^*)$ and the solution itself, based on the local Lipschitz property of $\mathcal{S}_{\mathrm{KKT}}$. Although the proof closely follows that of \cite[Theorem 3.1]{KS18}, we include it here for completeness.

\begin{theorem}\label{thm:pertkkt-resi}
	Let $(x^*, y^*, z^*)$ be a KKT point for problem (\ref{eq:prime}). The following statements are equivalent:
	\begin{itemize}
		\item [(i)] There exists a neighborhood $\mathcal{U}$ around $x^*$ and a constant $\kappa > 0$ such that for any perturbation $q = (\hat{a}, b, c) \in \mathbb{R}^n \times \mathbb{Y} \times \mathbb{Z}$ close to $(0, 0, 0)$, any solution $\left(x_q, y_q, z_q\right) \in \mathcal{U} \times \mathbb{Y} \times \mathbb{Z}$ of the perturbed KKT system (\ref{eq:man-kkt-perturb2}) satisfies:
		\[
		d(x_q, x^*) + \operatorname{dist}\left((y_q, z_q), M(x^*, 0, 0, 0)\right) \leq \kappa \|q\|_{\mathbb{R}^n \times \mathbb{Y} \times \mathbb{Z}}.
		\]
		\item [(ii)] There exists a neighborhood $\mathcal{U}$ around $x^*$ and a constant $\kappa > 0$ such that for any $(x, y, z) \in \mathcal{U} \times \mathbb{Y} \times \mathbb{Z}$ with $R(x, y, z)$ sufficiently small, the following inequality holds:
		\begin{equation}\label{eq:dist-resi}
		d(x,x^*)+\operatorname{dist}\left((y,z),M(x^*,0,0, 0)\right)  \leq \kappa R(x, y,z).
	\end{equation}
	\end{itemize}
\end{theorem}
\proof (ii) $\Rightarrow$ (i) Let $q=(\hat{a},b,c)\in\mathbb{R}^n\times\mathbb{Y}\times\mathbb{Z}$. Suppose that $(x_q, y_q,z_q)$ is the solution of (\ref{eq:man-kkt-perturb2}) for $q$, then
$$
\operatorname{grad} _xL(x_q, y_q,z_q)=(D\varphi(x_q))^{-1}G_{\varphi(x_q)}^{-1}\hat{a}.
$$
Since $(D\varphi(x))^{-1}$ and $G_{\varphi(x)}^{-1}$ are smooth functions of $x$, there exist $L>0$ and a neighborhood $\mathcal{ U }_{x^*}$ of $x^*$, such that for any $x\in \mathcal{ U }_{x^*}$, we have
$$
\|(D\varphi(x))^{-1}G_{\varphi(x)}^{-1}\hat{a}\|\leq L\|\hat{a}\|.
$$
By (\ref{eq:man-kkt-perturb2}), $g_1\left(x_q\right)+b=\operatorname{prox}_{\theta}\left(g_1\left(x_q\right)+b+y_q\right)$ and $g_2\left(x_q\right)+c=\Pi_{\mathcal{Q}}\left(g_2\left(x_q\right)+c+z_q\right)$. Thus, we obtain that
$$
\begin{aligned}
	&\left\|g_1\left(x_q\right)-\operatorname{prox}_{\theta}\left(g_1\left(x_q\right)+y_q\right)\right\|\\
	=&\left\|g_1\left(x_q\right)-\operatorname{prox}_{\theta}\left(g_1\left(x_q\right)+y_q\right)\right\| -\left\| g_1\left(x_q\right)+b-\operatorname{prox}_{\theta}\left(g_1\left(x_q\right)+b+y_q\right) \right\|\\
	\leq&\left\|g_1\left(x_q\right)-\operatorname{prox}_{\theta}\left(g_1\left(x_q\right)+y_q\right)-\left( g_1\left(x_q\right)+b-\operatorname{prox}_{\theta}\left(g_1\left(x_q\right)+b+y_q\right)\right) \right\|\\
	\leq&\left\|g_1(x_q)-(g_1(x_q)+b)\right\|=\left\|b\right\|,
\end{aligned}
$$
where the last inequality holds by Moreau decomposition \cite{M65} and the fact that the proximal operator is 1-Lipschitz (see \cite[Proposition 12.19]{RW09}). Similarly, we have $\left\|g_2\left(x_q\right)-\Pi_{\mathcal{Q}}\left(g_2\left(x_q\right)+z_q\right)\right\| \leq\|c\|$.
Therefore,
\begin{equation}\label{eq:resi-q}
	\begin{array}{lll}
		R\left(x_q, y_q,z_q\right)&=&\|\operatorname{grad} _xL(x_q, y_q,z_q)\|+\left\|g_1\left(x_q\right)-\operatorname{prox}_{\theta}\left(g_1\left(x_q\right)+y_q\right)\right\|\\
		&&+\left\|g_2\left(x_q\right)-\Pi_{\mathcal{Q}}\left(g_2\left(x_q\right)+z_q\right)\right\|\\
		&\leq& L\|\hat{a}\|+\|b\|+\|c\| \leq  (L+1)\|q\| .
	\end{array}
\end{equation}
Thus, taking $q$ sufficiently close to zero and substituting (\ref{eq:resi-q}) into (\ref{eq:dist-resi}), we now obtain
$$
d(x_q,x^*)+\operatorname{dist}\left((y_q,z_q), M(x^*, 0, 0,0)\right) \leq \kappa R\left(x_q,y_q,z_q\right) \leq \kappa(L+1)\|q\|.
$$
(i) $\Rightarrow$ (ii) Let $(x, y,z) \in \mathcal{U} \times \mathbb{Y}\times \mathbb{Z}$ and we define
$$
\tilde{g}_1:=\operatorname{prox}_{\theta}(g_1(x)+y), \quad \tilde{y}:=g_1(x)+y-\tilde{g}_1,\quad\tilde{g}_2:=\Pi_{\mathcal{Q}}(g_2(x)+z), \quad \tilde{z}:=g_2(x)+z-\tilde{g}_2.
$$
Now take $\hat{a}=G_{\varphi(x)}D\varphi(x)\operatorname{grad} _xL(x,\tilde{y},\tilde{z})$, $b=g_1(x)-\tilde{g}_1$ and $c=g_2(x)-\tilde{g}_2$. It is easy to obtain that $\tilde{y}\in\partial\theta(g_1(x)+b)$ and $\tilde{z} \in \mathcal{N}_{q}(g_2(x)+c)$. Hence, $(x,\tilde{y}, \tilde{z})$ is the solution of (\ref{eq:man-kkt-perturb2}) with $(\hat{a},b,c)$. We further have
$$
\|\tilde{y}-y\|=\|b\|=\|\tilde{g}_1-g_1(x)\|=\left\|g_1(x)-\operatorname{prox}_{\theta}(g_1(x)+y)\right\| \leq R(x, y,z)
$$
$$
\text{and}\quad\|\tilde{z}-z\|=\|c\|=\|\tilde{g}_2-g_2(x)\|=\left\|g_2(x)-\Pi_{\mathcal{Q}}(g_2(x)+z)\right\|\leq R(x, y,z)
$$
by the definition of $\tilde{y},\tilde{z}$ and $R(x, y,z)$. Now assume that $\left\|Dg_1(x)^{*}\right\| \leq c_{1}$ and $\left\|Dg_2(x)^{*}\right\| \leq c_{2}$ for all $x \in \mathcal{U}$ with some $c_{1},c_2\geq 0 $, we have
$$
\begin{array}{ll}
	\|(\hat{a},b,c)\|&=\|G_{\varphi(x)}D\varphi(x)\operatorname{grad} _xL(x,\tilde{y},\tilde{z})\|+\|b\|+\|c\|\\
	& \leq L_1\|\operatorname{grad} _xL(x,\tilde{y},\tilde{z})\|+\|b\|+\|c\|\\
	&\leq L_1\|\operatorname{grad} _xL(x,y,z)\|+(L_1 c_1+1)\|b\|+(L_1c_2+1)\|c\|\\  &\leq\left(L_1(1+c_{1}+c_2)+2\right) R(x, y,z).
\end{array}
$$
The first inequality is obtained from the smoothness of  $D\varphi(\cdot)$ and $G_{\varphi(\cdot)}$ in $\mathcal{U}$ which implies that there exists $L_1>0$, such that for any $x\in\mathcal{U}$, $
\| G_{\varphi(x)}D \varphi(x) \operatorname{grad}_x L(x,\tilde{y},\tilde{z}) \|\leq L_{1} \| \operatorname{grad}_x L(x,\tilde{y},\tilde{z})\|$. Therefore, taking $R(x,y,z)$ close enough to 0, we can then apply (i) to $(x,\tilde{y},\tilde{z})$ and obtain
$$
d(x,x^*)+\operatorname{dist}\left((\tilde{y},\tilde{z}), M(x^*, 0,0, 0)\right) \leq \kappa\|(\hat{a},b,c)\| \leq \kappa\left(L_1(1+c_{1}+c_2)+2\right) R(x, y,z) .
$$
Moreover, by $\|\tilde{y}-y\|\leq R(x, y,z)$ and $\|\tilde{z}-z\| \leq R(x, y,z)$, 
$$
\operatorname{dist}\left((y,z), M(x^*,0, 0, 0)\right)\leq \operatorname{dist}\left((\tilde{y},\tilde{z}),M(x^*,0, 0, 0)\right) + 2R(x, y,z).
$$  
We finally obtain
$$
d(x,x^*)+\operatorname{dist}\left((y,z), M(x^*,0,0, 0)\right) \leq [\kappa\left(L_1(1+c_{1}+c_2)+2\right) +2] R(x, y,z).
$$
The proof is then complete.\qed

For any $\sigma\geq 0,$ we denote by $M_\sigma(x)$ the set of $(y,z)$ satisfying the following relationships: 
\begin{equation}\label{eq:perturb-multi}
	\|\operatorname{grad} _xL(x,y,z)\|\leq \sigma,\quad y\in\partial\theta(g_1(x)),\quad z\in\mathcal{ N }_{\mathcal{Q}}(g_2(x)).
\end{equation}
The following proposition is a direct conclusion from \cite[Theorem 4.43]{BS13}.

\begin{proposition}\label{pro:mult-bound}
	Suppose that M-RCQ holds at $x^*$. Then for any $\sigma\geq0$, the solution sets $M_\sigma(x)$  are 
	bounded for all $x$ in a neighborhood of $x^*$.
\end{proposition}
\proof Let us denote $\widehat{M}_\sigma(\hat{x})$ as the set of $(y,z)$ satisfying 
\begin{equation}\label{eq:perturb-multi2}
	\|\nabla_{\hat{x}}\widehat{L}(\hat{x},y,z)\|\leq \sigma ,\quad y\in\partial\theta(\hat{g}_1(\hat{x})),\quad z\in\mathcal{ N }_{\mathcal{Q}}(\hat{g}_2(\hat{x})),
\end{equation}
where $\hat{x}=\varphi(x)$ and $\widehat{L}=L\circ\varphi^{-1}$. By (\ref{eq:co-ex}) and the Lipschitz property of $D\varphi(x)$ and $G_{\varphi(x)}$ at $x^*$, if $(y,z)$ satisfying (\ref{eq:perturb-multi}) for some $\sigma$, then there exists $\alpha\geq 0$, such that for any $x$ in a neighborhood $\mathcal{ U }_{x^*}$,  $(y,z)\in\widehat{M}_{\alpha\sigma}(\hat{x})$. Thus, $M_{\sigma}(x)\subseteq \widehat{M}_{\alpha\sigma}(\hat{x}).$ Using \cite[Theorem 4.43]{BS13}, we have the boundness of $\widehat{M}_{\alpha\sigma}(\hat{x})$, which implies the boundness of $M_\sigma(x)$.\qed

With this proposition, we can establish the following local error bound around the KKT point $(x^*,y^*,z^*)$ in terms of the KKT residual mapping. 

\begin{theorem}\label{thm:errorbound}
	Let $(x^*, y^*, z^*)$ be a KKT point of problem (\ref{eq:man-perturb}) with $(\hat{a}, b, c) = (0, 0, 0)$ that satisfies both the M-SOSC (\ref{eq:man-sosc-perturb}) and M-SRCQ (\ref{eq:man-srcq-perturb}) conditions. Then $M(x^*, 0, 0, 0) = \{(y^*, z^*)\}$. Furthermore, for any $(x, y, z) \in \mathcal{M} \times \mathbb{Y} \times \mathbb{Z}$ with $x$ sufficiently close to $x^*$ and $R(x, y, z)$ sufficiently small, there exist constants $c_1, c_2 > 0$ such that:
	\begin{equation}\label{eq:2errorbound}
	c_{1} R(x,y,z) \leq	d(x,x^*)+\|y-y^*\|+\|z-z^*\| \leq c_2 R(x,y,z).
\end{equation}
\end{theorem}
\proof Let $\{(\hat{a}^k, b^k, c^k)\}$ be a sequence converging to 0, and let $\{(x^k, y^k, z^k)\}$ be a sequence of solutions to (\ref{eq:man-kkt-perturb2}) corresponding to $\{(\hat{a}^k, b^k, c^k)\}$, such that $x^k \to x^*$. By Proposition \ref{pro:mult-bound}, $\{(y^k, z^k)\}$ is bounded, which implies that every accumulation point of $\{(y^k, z^k)\}$ is the multiplier corresponding to $x^*$. Since M-SRCQ holds at $x^*$ with respect to $(y^*, z^*)$, Theorem \ref{thm:man-srcq} shows that $M(x^*) = \{(y^*, z^*)\}$, thus $(y^k, z^k)$ converges to $(y^*, z^*)$. By Remark \ref{rm:iso-equi-naturalmapping} and the definition of $R(x^k, y^k, z^k)$, there exists a constant $\kappa > 0$ such that for sufficiently large $k$,
$$
d(x^k, x^*) + \|y^k - y^*\| + \|z^k - z^*\| \leq \kappa R(x^k, y^k, z^k) \leq \kappa \|(\hat{a}^k, b^k, c^k)\|.
$$
The second inequality follows from Theorem \ref{thm:pertkkt-resi} for some $c_2 > 0$.

Note that the function $R$ is locally Lipschitz continuous with respect to $x$ and globally with respect to $(y,z)$. Therefore, 
the left term holds for a suitable constant $c_1>0$.\qed

We derive the following lemma to analyze the behavior of local minimizers of the augmented Lagrangian function in a neighborhood of $ x^* $. This lemma extends the result from \cite[Lemma 4.1]{KS19} from Euclidean to Riemannian settings.

\begin{lemma}\label{lem:exist-mini}
	Let $(x^*, y^*, z^*)$ be a KKT point of problem (\ref{eq:man-perturb}) with $(\hat{a}, b, c) = (0, 0, 0)$ that satisfies the M-SOSC (\ref{eq:man-sosc-perturb}) condition, and let $B \subseteq \mathbb{Y} \times \mathbb{Z}$ be a bounded set. Then there exist constants $\bar{\rho} > 0$ and $r > 0$ such that for every $\rho \geq \bar{\rho}$ and every $(w, p) \in B$, the function $L_{\rho}(x, w, p)$ has a local minimizer $x = x_{\rho}(w, p)$ within the ball $B_r(x^*)$ (where $B_r(x^*)$ denotes the closed ball of radius $r$ around $x^*$ with respect to the Riemannian distance defined in \eqref{eq:def-R-distance} on $\mathcal{M}$). Furthermore, $x_{\rho} $  uniformly converges to $x^*$ on $B$ as $\rho \to \infty$.
\end{lemma}
\proof  For each $\rho > 0$ and $(w, p) \in B$, by the compactness of $B_{r}(x^*)$, there exists $x = x_{\rho}(w, p)$ which is the solution of the problem
$$
\min _{x}\; L_{\rho}(x,w,p) \quad\text {s.t.}\; x \in B_{r}(x^*).
$$
Given that the M-SOSC holds, there exists an $r > 0$ such that $x^*$ is a strict local solution of this problem within $B_{r}(x^*)$. Assume, for contradiction, that $x_{\rho}$ does not uniformly converge to $x^*$. Then there exist $\epsilon > 0$, a sequence $\rho^{k} \to \infty$, and a sequence $\{(w^k, p^k)\} \subseteq B$ such that $d(x_{\rho^{k}}(w^k,p^k),x^*)>\epsilon$ for all $k$.
Since $B_{r}(x^*)$ is compact, the sequence $\{x_{\rho^k}(w^k, p^k)\}$ has an accumulation point $\tilde{x}$. By the properties of $g_2(x^*) \in \mathcal{Q}$, we have for all $k$,
$$
\theta^{\rho^k}\left(g_1(x^*)+\dfrac{w^k}{\rho^k}\right)\leq \theta(g_1(x^*))+\dfrac{\|w^k\|^2}{2\rho^{k}}\quad\text{and}\quad
\dfrac{\rho^{k}}{2}\operatorname{dist}^2\left(g_2(x^*)+\dfrac{p^k}{\rho^{k}},\mathcal{Q}\right)\leq\dfrac{\|p^k\|^2}{2\rho^{k}}.
$$
Thus, we have
\begin{eqnarray}
	&&f\left(x^{k}\right)+\theta^{\rho^k}\left(g_1(x^k)+\dfrac{w^k}{\rho^k}\right)+\dfrac{\rho^{k}}{2}\operatorname{dist}^2(g_2(x^k)+\dfrac{p^k}{\rho^k},\mathcal{Q}) \nonumber \\
	&\leq & L_{\rho^{k}}\left(x^*, w^{k},p^k\right)\leq  f(x^*)+\theta(g_1(x^*))+\dfrac{\left\|w^{k}\right\|^{2}}{2\rho^{k}}+\dfrac{\|p^k\|^2 }{2\rho^{k}}.\label{eq:mini-xk}
\end{eqnarray}
Since $$\theta^{\rho^k}\left(g_1(x^k)+\dfrac{w^k}{\rho^k}\right)=\theta\left(\operatorname{prox}_{\theta/\rho^{k}}(g_1(x^k)+\dfrac{w^k}{\rho^k})\right)+\dfrac{\rho^k}{2}\|g_1(x^k)+\dfrac{w^k}{\rho^k}-\operatorname{prox}_{\theta/\rho^{k}}(g_1(x^k)+\dfrac{w^k}{\rho^k})\|,$$
by taking $\rho^{k}\to \infty$, we know from (\ref{eq:mini-xk})  that
$$\|g_1(x^k)+\dfrac{w^k}{\rho^k}-\operatorname{prox}_{\theta/\rho^{k}}(g_1(x^k)+\dfrac{w^k}{\rho^k})\|\to 0\quad  \mbox{and}\quad \operatorname{dist}(g_2(x^k)+\dfrac{p^k}{\rho^k},\mathcal{Q})\to 0.$$ Hence, we have $g_1(\tilde{x})$ is an auccumulation point of  $\operatorname{prox}_{\theta/\rho^{k}}(g_1(x^k)+\dfrac{w^k}{\rho^k})$ and $g_2(\tilde{x})\in \mathcal{Q}$. Moreover, (\ref{eq:mini-xk}) also yields 
$\limsup _{k \rightarrow \infty} f\left(x^{k}\right)+\theta\left( \operatorname{prox}_{\theta/\rho^{k}}(g_1(x^k)+\dfrac{w^k}{\rho^k})\right)  \leq f(x^*)+\theta(g_1(x^*))$. 
Therefore, $f(\tilde{x})+\theta(g_1(\tilde{x}))\leq f(x^*)+\theta(g_1(x^*))$, which means that $\tilde{x}=x^*$ since $x^*$ is the strict solution in $B_{r}(x^*)$. This contradicts the assumption, and the proof is thus complete.\qed

\begin{remark}\label{rm:inexact-solution}
	The uniform convergence implies the existence of a $\bar{\rho}>0$ such that $x_{\rho}(w,p)$ lie in the interior of $ B_{r}(x^*)$, $(w,p)\in B$. Now taking an appropriate small $\epsilon_{\rho}>0$, by the continuity of $L_{\rho}(x, w,p)$ and $\operatorname{grad}_{x}L_{\rho}(x, w,p)$ for $x$, inequality (\ref{eq:subproblem}) and $L_{\rho}(\tilde{x}_{\rho},w,p)\leq L_{\rho}(x_{\rho},w,p)+\epsilon_{\rho}$ hold simultaneously for any $\tilde{x}_{\rho}$ sufficiently close to $x_{\rho}$. Then, we obtain that $\tilde{x}_{\rho}\to x^*$ if $\epsilon_\rho\to 0$ as $\rho\to \infty$.
\end{remark}

Next, we establish the local convergence result for the inexact Riemannian ALM (Algorithm \ref{alg:alm}), which is inspired by \cite[Theorem 4.2]{KS19}. The proof can be applied to the linear convergence of exact ALM when we take each $\epsilon_k=0$.

\begin{theorem}\label{thm:converge}
	Let $(x^*, y^*, z^*)$ be a KKT point of problem (\ref{eq:man-perturb}) with $(\hat{a}, b, c) = (0, 0, 0)$ that satisfies both the M-SOSC (\ref{eq:man-sosc-perturb}) and M-SRCQ (\ref{eq:man-srcq-perturb}) conditions. Suppose $x^{k+1}$ is an approximate minimizer as described in Remark \ref{rm:inexact-solution} for sufficiently large $k$. Then there exists a $\bar{\rho} > 0$ such that if $\rho^k \geq \bar{\rho}$ for sufficiently large $k$, it follows that $\left(x^k, y^k, z^k\right) \to (x^*, y^*, z^*)$. Additionally, if $\epsilon_k = o(R(x^k, y^k, z^k))$ and $(w^k, p^k) = (y^k, z^k)$ for sufficiently large $k$, then there exists a constant $c > 0$, such that
	\begin{equation}\label{eq:conver-rate}
	d(x^{k+1},x^*)+\|y^{k+1}-y^*\|+\|z^{k+1}-z^*\|\leq \dfrac{c}{\rho^{k}}\left(d(x^k,x^*)+\|y^{k}-y^*\|+\|z^{k}-z^*\|\right)
\end{equation}
	holds for all sufficiently large $k$. Furthermore, the sequence $\{\rho^k\}$ remains bounded.
\end{theorem}
\proof By choosing $\bar{\rho} > 0$ sufficiently large, the selected $x^k$ can be made sufficiently close to $x^*$. According to the two-sided error bound given in Theorem \ref{thm:errorbound}, for any $\rho^k \geq \bar{\rho}$, $x^k$ lies in a neighborhood of $x^*$ where the error bound given by (\ref{eq:2errorbound}) is valid. Consequently, the convergence of $R_{k} := R(x^k, y^k, z^k) \to 0$ implies that $(x^k, y^k, z^k) \to (x^*, y^*, z^*)$.
From the definition of $(x^{k}, y^{k},z^k)$, we know that
$$
\operatorname{grad} _xL(x^{k},y^{k},z^k)=\operatorname{grad} _x L _{\rho^k}(x^k,w^{k},p^{k}).
$$
Now let $
p^{k+1}:=\operatorname{prox}_{\theta/\rho^k}\left(g_1\left(x^{k+1}\right)+w^{k}/ \rho^{k}\right)\text{ and }
s^{k+1}:=\Pi_{\mathcal{Q}}\left(g_2\left(x^{k+1}\right)+p^{k} / \rho^{k}\right)$. 
Then $y^{k+1} \in \partial\theta(p^{k+1})$ and $z^{k+1}\in\mathcal{ N }_{\mathcal{Q}}(s^{k+1})$, which imply that $p^{k+1}=\operatorname{prox}_{\theta/\rho^k}\left(p^{k+1}+y^{k+1}\right) $ and $s^{k+1}=\Pi_{\mathcal{Q}}\left(s^{k+1}+z^{k+1}\right) $. Therefore, we have
\begin{eqnarray}
	&&\left\|g_1\left(x^{k+1}\right)-\operatorname{prox}_{\theta/\rho^k}(g_1(x^{k+1})+y^{k+1})\right\|\nonumber\\
	&=&\left\|g_1\left(x^{k+1}\right)-\operatorname{prox}_{\theta/\rho^k}(g_1(x^{k+1})+y^{k+1})\right\|-\left\|p^{k+1}-\operatorname{prox}_{\theta/\rho^k}\left(p^{k+1}+y^{k+1}\right)\right\| \nonumber\\
	&\leq &\left\|g_1\left(x^{k+1}\right)-\operatorname{prox}_{\theta/\rho^k}(g_1(x^{k+1})+y^{k+1})-p^{k+1}+\operatorname{prox}_{\theta/\rho^k}\left(p^{k+1}+y^{k+1}\right)  \right\| \nonumber\\
	&\leq&\left\|g_1\left(x^{k+1}\right)-p^{k+1}\right\|.\label{eq:g1-prox}
\end{eqnarray}
The last inequality is obtained by using Moreau decomposition and the Lipschitz property of the proximal operator. Similarly, we have 
\begin{equation}\label{eq:g2-pi}
	\left\|g_2\left(x^{k+1}\right)-\Pi_{\mathcal{Q}}(g_2(x^{k+1})+z^{k+1})\right\|\leq\left\|g_2\left(x^{k+1}\right)-s^{k+1}\right\|.
\end{equation}

When $\{\rho^k\}$ is bounded, then by  (\ref{eq:update-rho}),
$\left\|g_1\left(x^{k+1}\right)-p^{k+1}\right\|$ and $
\left\|g_2\left(x^{k+1}\right)-s^{k+1}\right\|$ both converge to zero.
If $\rho^{k} \rightarrow \infty$,  by Lemma \ref{lem:exist-mini} and Remark \ref{rm:inexact-solution}, we have $x^{k+1} \rightarrow x^*$,
$$
\left\|g_1\left(x^{k+1}\right)-p^{k+1}\right\| =\left\|g_1(x^{k+1})-\operatorname{prox}_{\theta/\rho^k}\left(g_1(x^{k+1})+\dfrac{w^k}{\rho^k}\right)\right\| \rightarrow 0 
$$
and
$$
\left\|g_2\left(x^{k+1}\right)-s^{k+1}\right\| \leq\left\|s^{k+1}-\Pi_{\mathcal{Q}}\left(g_2\left(x^{k+1}\right)\right)\right\|+\operatorname{dist}\left(g_2\left(x^{k+1}\right),\mathcal{Q}\right) \rightarrow 0.
$$
Therefore, we know that
\begin{equation}\label{eq:resi-to-0}
	R_{k+1} \leq \epsilon_{k+1}+\left\|g_1\left(x^{k+1}\right)-\operatorname{prox}_{\theta/\rho^k}(g_1(x^{k+1})+y^{k+1})\right\|+	\left\|g_2\left(x^{k+1}\right)-\Pi_{\mathcal{Q}}(g_2(x^{k+1})+z^{k+1})\right\|\rightarrow 0,
\end{equation}
which means that $\left(x^{k},y^{k},z^k\right) \rightarrow(x^*,y^*,z^*)$.

Now we are considering the convergence rate. Assume that for sufficiently large $k$, $(w^{k},p^k)=(y^{k},z^k)$ and $\epsilon_k=o(R(x^{k}, y^{k},z^k))$ . Using (\ref{eq:resi-to-0}) and the definition of $(y^{k+1},z^{k+1})$, we have
\begin{eqnarray}
	R_{k+1} &\leq& \epsilon_{k+1}+\left\|g_1\left(x^{k+1}\right)-p^{k+1}\right\|+\left\|g_2\left(x^{k+1}\right)-s^{k+1}\right\| \nonumber \\
	&=&\epsilon_{k+1}+\dfrac{\left\|y^{k+1}-y^{k}\right\|}{\rho^{k}}+\dfrac{\left\|z^{k+1}-z^k\right\|}{\rho^{k}} \nonumber\\
	&\leq& \epsilon_{k+1}+\dfrac{1}{\rho^{k}}\left(\left\|y^{k+1}-y^*\right\|+\left\|y^{k}-y^*\right\|+\left\|z^{k+1}-z^*\right\|+\left\|z^{k}-z^*\right\|\right).\label{eq:resi-multi}
\end{eqnarray}
(\ref{eq:2errorbound}) implies that there exists $c_{1}>0$, such that $\left\|y^{k}-y^*\|+\|z^k-z^*\right\| \leq c_{1} R_{k}$.  Since $\epsilon_k=o(R(x^k,p^k,y^{k},z^k))$, there exist $\alpha<\dfrac{c_1}{3\bar{\rho}}$ and $\delta>0$, such that 
$$
\epsilon_{k}\leq \alpha R(x^k,y^{k},z^k) \quad\forall\,(x^k,y^{k},z^k)\text{ satisfying }d((x^k,y^{k},z^k),(x^*,y^*,z^*))\leq \delta.
$$
It follows that $R_{k+1} \leq\dfrac{c_1}{\rho^k}\left(R_{k+1}+R_{k}\right)+\alpha R_{k+1}$. Now increase $\bar{\rho}$ until $1-\dfrac{c_1}{\rho^k}-\alpha>1 / 2$, then
$
\left(1-\dfrac{c_{1}}{\rho^{k}}-\alpha \right) R_{k+1} \leq \dfrac{c_{1}}{\rho^{k}} R_{k}
$
implies that $R_{k+1}<\left(2 c_{1} / \rho^{k}\right) R_{k}$. The linear convergence is now obtained by (\ref{eq:2errorbound}).

Let $V_{k+1}:=V\left( x^{k+1}, y^{k+1},z^{k+1}, \rho^k\right) $. The boundness of $\left\{\rho^{k}\right\}$ can be obtained if $V_{k+1} \leq \tau V_{k}$ holds for $k$ sufficiently large. By the equality in (\ref{eq:resi-multi}), we have $2V_{k+1} \geq R_{k+1}-\alpha R_{k}$ and  
\begin{eqnarray*}
	V_{k+1}&=&\max\left(\dfrac{\left\|y^{k+1}-y^{k}\right\|}{\rho^{k}},\dfrac{\left\|z^{k+1}-z^k\right\|}{\rho^{k}}\right) \\
	&\leq& \dfrac{1}{\rho^{k}}\left(\left\|y^{k+1}-y^*\right\|+\left\|y^{k}-y^*\right\|+\left\|z^{k+1}-z^*\right\|+\left\|z^{k}-z^*\right\|\right)\leq \dfrac{c_{1}}{\rho^{k}}\left(R_{k+1}+R_{k}\right) .
\end{eqnarray*}
Thus, these inequalities yield
$$
\dfrac{V_{k+1}}{V_{k}} \leq \dfrac{2c_{1}\left(R_{k+1}+R_{k}\right) }{\rho^{k} \left( R_{k}-\alpha R_{k-1}\right) }\leq \dfrac{2c_{1}}{\rho^{k}}\left(1+\dfrac{(2c_1/\rho^{k})R_{k}+\alpha R_{k-1}}{R_{k}-\alpha R_{k-1}}\right).
$$
When $\rho^{k} \rightarrow \infty$, it can be seen that $V_{k+1} / V_{k} \rightarrow 0$, which implies that $\{\rho^{k}\}$ is bounded.\qed

\begin{remark}\label{rm:proof-wp-bounded}
If $\{(y^k, z^k)\}$ is bounded, then we can choose $B$ as a bounded set that contains $\{(y^k, z^k)\}$ for all $k$. This choice ensures that the condition in Theorem \ref{thm:converge}, which requires $(w^k, p^k)$ to be equal to $(y^k, z^k)$ for sufficiently large $k$, is satisfied. Consequently, we can simply select $(w^k, p^k)$ to be $(y^k, z^k)$.
\end{remark}

\begin{remark}\label{rm:penalty-rate}
	The convergence rates given in (\ref{eq:conver-rate}) depend on the constant $c$ and the penalty parameters $\rho^{k}$, where $c$ is determined by the problem and $\rho^{k}$ can be chosen dynamically. It is shown that the fast linear rate can be achieved if we increase the penalty parameters $\rho^{k}$, and the rates become asymptotically superlinear if $\rho^{k}\to\infty$.
\end{remark}

\begin{remark}\label{rm:example-sosc-srcq-alm}
	The well behavior of  Riemannian ALM  solving (\ref{exam:sphere}) can now be explained, since the  M-SOSC condition is satisfied at $x^*$ and M-SRCQ holds at $x^*$ with respect to $(y^*,z^*)$. According to Theorem \ref{thm:converge}, these conditions ensure that the iteration sequence converges linearly for sufficiently large $k$.
\end{remark}

\section{Applications and numerical experiments}\label{sec:numerical}
It is established by the two-sided error bound (\ref{eq:2errorbound}) that the linear convergence of the KKT residual $R(x,y,z)$ can imply the linear convergence of the iteration residual. We will use 
 $R(x,y,z)$  to illustrate the convergence rate for two reasons. First, it is challenging to obtain the local solution and corresponding multipliers for randomly generated cases. Second, the distance functions for most manifolds are unknown, making it unconvincing to rely solely on Euclidean distance.
In the following experiments, the ALM penalty parameters $\rho^{k}$ are tuned adaptively to guarantee smooth implementation.
%
All codes are implemented in Matlab (R2021b) and all the numerical experiments are run under a  64-bit macOS on an Intel Core i5 2.4GHz CPU with 16GB memory.

\subsection{Nonsmooth optimization on sphere}
Set $g_1(x)=x$,  $\theta(\cdot)=\mu\|\cdot\|_1$ and $\mathcal{M}$ to be the sphere ${\cal S}^{n-1}:=\{x\in\mathbb{R}^n\mid x^Tx=1\}$. Consider the following nonsmooth optimizations on the unit sphere
\begin{equation}\label{eq:sphere}
	\min_{x \in \mathcal{M}=\mathcal{S}^{n-1}}  f(x)+\mu\|x\|_{1},
\end{equation}
where $f:\mathcal{S}^{n-1}\rightarrow \mathbb{R}$ is a smooth function. It is well-known that tangent space and normal space for $\mathcal{S}^{n-1}$ at a point $x$ are given by $T_{x}\mathcal{ M }=\left\lbrace \xi\mid x^{\top}\xi=0 \right\rbrace$ and $N_{x}\mathcal{M}=\left\lbrace ax\mid a\in\mathbb{R}\right\rbrace $,  respectively. The projection onto the tangent space is 
$\Pi_{x}(z)=z-xx^{\top}z$ and $\operatorname{grad}f(x)=\Pi_{x}(\nabla f(x))$.  Follows from $g_1(x)=x$, $Dg_1(x)T_{x}\mathcal{M}=T_{x}\mathcal{M}$ and $Dg_1(x)^{*}y=\Pi_{x}(y)$ for any $y\in\mathbb{R}^n$. Since $\operatorname{dom}\theta=\mathbb{R}^n$, it can be checked directly that the M-RCQ condition (\ref{eq:man-rcq-perturb}) with respect to \eqref{eq:sphere} is satisfied at any feasible point of (\ref{eq:sphere}).\footnote{The corresponding M-RCQ condition (\ref{eq:man-rcq-perturb}) also holds at any feasible point of (\ref{eq:sphere}) where the manifold $\mathcal{ M }$ in (\ref{eq:sphere}) is replaced by the Stiefel manifold $\operatorname{St}(n,p)=\{X\in \mathbb{R}^{n\times p}\mid X^{\top}X=I_p\}$, e.g., the sparse principal component analysis (SPCA) problem (see \cite{JTU03,ZHT06}).}  
Furthermore, the M-SRCQ condition (\ref{eq:man-srcq-perturb}) with respect to \eqref{eq:sphere} is said to hold at a stationary point $x$ with a multiplier $y$ if 
\begin{equation}\label{eq:sphere-srcq}
	T_{x}\mathcal{ M}+	\mathcal{C}_{\theta,g_1}(x, y)=\mathbb{R}^{n},
\end{equation}
where  
\begin{equation}\label{eq:sphere-critical-theta}
	\mathcal{C}_{\theta,g_1}(x,y )=\left\{d \in \mathbb{R}^{n} \mid \theta^{\downarrow}(x ; d)=\langle d, y\rangle\right\rbrace \text{ with } y \text{ satisfying } \Pi_{x}(y) = -\Pi_{x}(\nabla f(x)).
\end{equation}
Since $\theta(\cdot)=\mu\|\cdot\|_1$, we know that  
$$\theta^{\downarrow}(x ; d)=\mu \sum_{x_{i}=0}|d_{i}|+\mu\sum_{x_{i}>0}d_{i}-\mu\sum_{x_{i}<0}d_{i}.  $$
Since the M-RCQ condition holds at any feasible point of \eqref{eq:sphere}, it follows from Theorem \ref{thm:man-rcq} that if $x^*$ is a local optimal solution of \eqref{eq:sphere}, then $x^*$ must be a stationary point, i.e., there exists $y\in\mathbb{R}^n$ such that
\begin{equation}\label{eq:sphere-kkt}
	\left\{\begin{array}{l}
		\Pi_{x^*}( \nabla f(x^*)+y )=0, \\
		y \in \mu\partial\|x^*\|_{1}.
	\end{array}\right.
\end{equation}
Moreover, we know from the following proposition that the M-SRCQ condition (\ref{eq:man-srcq-perturb}) for problem (\ref{eq:sphere}) is indeed satisfied at any KKT pair $(x^*,y^*)$, which implies that the multiplier is unique.

\begin{proposition}\label{pro:shpere-srcq}
    Let $x^*$ be a local optimal solution of (\ref{eq:sphere}) and $y^*$ be a corresponding multiplier satisfying (\ref{eq:sphere-kkt}). Then, the M-SRCQ condition \eqref{eq:sphere-srcq} holds for (\ref{eq:sphere}) at $x^*$ with respect to $y^*$, which implies the corresponding multiplier $y^*$ is unique.
\end{proposition}
\proof In order to show (\ref{eq:sphere-srcq}) holds, we only need to verify that the normal space $N_{x^*}\mathcal{M}$ at $x^*$ satisfies $N_{x^*}\mathcal{M}\subseteq \mathcal{ C}_{\theta,g_1}(x^*,y^*)$. Since $(x^*,y^*)$ satisfying (\ref{eq:sphere-kkt}), there exists $a$ such that $y^*+\nabla f(x^*)=ax^*$ and $y^*\in \mu \partial \|x^*\|_1$. Since $x^*$ is orthogonal, it has at least one nonzero component, and we denote by $x^*_{i}$. Then $a$ can be determined by $a=\dfrac{y^*_{i}+\nabla f(x^*)_{i}}{x^*_{i}}=\dfrac{\operatorname{sgn}(x^*_{i})\mu+\nabla f(x^*)_{i}}{x^*_{i}}$. Therefore, $y^*$ is uniquely defined by $a$. Taking any $d\in N_{x^*}\mathcal{M}$, there exists $\tilde{a}$ such that $d=\tilde{a}x^*$. Moreover, $d_i=0$ if $x^*_i=0$. Hence, it can be easily obtained that $\left\langle y^*,d \right\rangle =\sum_{i}y^*_id_i=\sum_{x^*_i>0}\mu d_i -\sum_{x^*_i<0}\mu d_i=\theta^{\downarrow}(x^*,d)$, which implies that $N_{x^*}\mathcal{M}\subseteq\mathcal{ C}_{\theta,g_1}(x^*,y^*)$. Thus, the M-SRCQ condition is satisfied at any KKT pair for (\ref{eq:sphere}). The uniqueness then follows from Theorem \ref{thm:man-srcq}, directly. \qed

It is clear that $\operatorname{epi}\theta$ with $\theta=
\mu\|\cdot\|_1$ is polyhedral. Thus, $\theta$ is $C^2$-cone reducible, then the M-SOSC condition (\ref{eq:man-sosc-perturb}) can be applied to problem (\ref{eq:sphere}). The critical cone at a stationary point $x$ is 
$$
\begin{aligned}
	\mathcal{C}(x)&=\left\lbrace \xi\in T_{x}\mathcal{M}\mid Df(x)\xi+(\theta\circ g_1)^{\circ}(x;\xi)=0 \right\rbrace\\
	&= \left\lbrace\xi\in T_{x}\mathcal{M}\mid \left\langle \nabla f(x),\xi\right\rangle +\theta^{\downarrow}(x;\xi)=0\right\rbrace.
\end{aligned}
$$
It is easy to compute that for any $\xi\in T_{x}\mathcal{M}$,
$$
\begin{aligned}
	\theta_{-}^{\downarrow \downarrow}\left(g_{1}(x) ; D g_{1}(x) \xi, w \right) & = \liminf_{t \downarrow 0 \atop w^{\prime} \rightarrow w} \dfrac{\mu \|x+t\xi+\frac{1}{2}t^2 w'\|_1-\mu\|x\|_1-t\theta^{\downarrow}(x,\xi)}{\frac{1}{2}t^2}\\
	& = \mu \sum_{x_{i}=0, \xi_{i}=0}|w_{i}|+\mu\sum_{x_{i}>0 \text{ or}\atop x_{i}=0,\xi_{i}>0}w_{i}-\mu\sum_{x_{i}<0\text{ or}\atop x_{i}=0,\xi_{i}<0}w_{i}.
\end{aligned}
$$
Therefore, for any multiplier $y$ at $x$,
$$
\begin{aligned}
	&\psi_{\left(g_{1}(x), D g_{1}(x) \xi\right)}^{*}(y) = \sup_{w} \left\lbrace \langle y,w\rangle-	\theta_{-}^{\downarrow \downarrow}\left(g_{1}(x) ; D g_{1}(x) \xi, w \right) \right\rbrace\\
	=& \sup_{w} \big\lbrace \sum_{x_{i}=0, \xi_{i}=0}y_{i}w_{i}-\mu|w_{i}|+\sum_{x_{i}>0 \text{ or}\atop x_{i}=0,\xi_{i}>0}y_{i}w_{i}-\mu w_{i}+\sum_{x_{i}<0\text{ or}\atop x_{i}=0,\xi_{i}<0}y_{i}w_{i}+\mu w_{i} \big\rbrace.
\end{aligned}
$$
Thus, the M-SOSC condition (\ref{eq:man-sosc-perturb}) for problem (\ref{eq:sphere}) at a stationary point $x^*$ with the unique multiplier $y^*$, takes the following form 
\begin{equation}\label{eq:sphere-sosc}
	\left\langle \xi,\operatorname{Hess}f(x^*)\xi \right\rangle -\psi_{\left(g_{1}(x^*), D g_{1}(x^*) \xi\right)}^{*}(y^*)>0 \quad \forall\,\xi\in \mathcal{C}(x^*)\backslash\{0\}.
\end{equation}

Unlike the M-SRCQ condition (\ref{eq:sphere-srcq}), the M-SOSC (\ref{eq:sphere-sosc}) may not hold for some cases. Nevertheless, we could still find some examples that fulfill this condition. Let us consider the following example of problem (\ref{eq:sphere}). Set $\mu=0.25$ and  $$f(x)=-\operatorname{tr}(x^{\top}A^{\top}Ax) \quad \mbox{with} \quad
A=\begin{bmatrix}
	10&     0  &    0  &  0  &    0\\
	0 & 25 &    0    &  0  &   0\\
	0  &   0  & 1.028&1.104 &   0\\
	0 &    0 &1.104&1.672  &     0\\
	0    &  0  &    0 &    0 & 8
\end{bmatrix}.$$
Then, $x^*=\left[ 
0,-1,0,0,0\right]^{\top} $ is an optimal solution of problem (\ref{eq:sphere}) with the multiplier $y^*=\left[ 
0,-\mu,0,0,0\right]^{\top}$. The M-SOSC condition is trivially satisfied since the critical cone $\mathcal{C}(x^*)=\{0\}$. 

\begin{figure}[ht]
	\centering
	\includegraphics[scale=0.5]{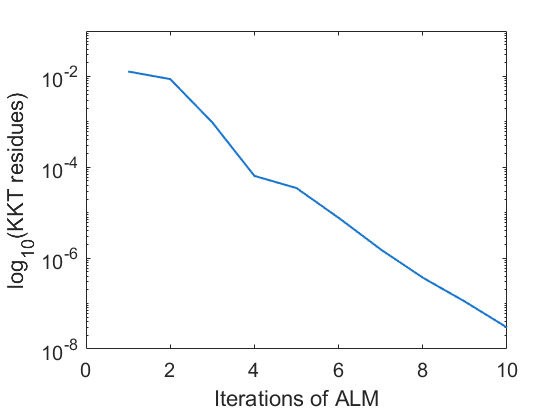}
	\caption{the KKT residues of problem (\ref{eq:sphere}) with randomly generated $A\in \mathbb{R}^{20\times 20}$.}
	\label{fig:SPCA}
\end{figure}

The detailed implementation\footnote{Matlab code is available at \url{https://github.com/miskcoo/almssn}.} and additional numerical results for the Riemannian ALM (Algorithm \ref{alg:alm}) are provided in \cite{ZBDZ21}. In this discussion, we use a simple example where the matrix $A \in \mathbb{R}^{20 \times 20}$ is randomly generated to illustrate the numerical performance of the Riemannian ALM. It is important to note that, since the exact solution for this random example is unknown, verifying the M-SOSC condition is challenging. Nevertheless, Figure \ref{fig:SPCA} demonstrates that the KKT residues converge linearly, even though the M-SOSC condition may not be satisfied.

\subsection{Robust matrix completion}\label{subsec:rmc}
Consider the robust matrix completion problem proposed in Example \ref{exam:rmc}.  For a given $A\in \mathbb{R}^{m\times n}$, let $g_1(X)=P_{\Omega}(X-A)$ and $\theta(\cdot)=\mu\|\cdot\|_{1}$. Here,  $P_{\Omega}$ is the projector defined by $\left(P_{\Omega}(X)\right)_{i j}=X_{i j}$ if $(i, j) \in \Omega$ and $0$ otherwise. By setting $\mathcal{M}=Fr(m,n,r):=\{X\in \mathbb{R}^{m\times n}:\operatorname{rank}(X)=r\}$, we obtain the following robust matrix completion problem
\begin{equation}\label{eq:rmc-full}
	\begin{array}{ll}
		\min _{X \in \mathbb{R}^{m \times n}} & \left\|P_{\Omega}(X-A)\right\|_{1} \\
		\text { s.t. } & X \in Fr(m,n,r).
	\end{array}
\end{equation}
Compared with the matrix completion using the Frobenius norm as an objective function, the $l_1$-norm is expected to due with the inexact data $A$ with some extreme outliers. 

It is known in \cite{V13} that the tangent space of $\mathcal{M}=Fr(m,n,r)$ at a point $X=USV^{\top}$ is
$$
\begin{aligned}
	T_{X} \mathcal{M}&=\left\{\left[\begin{array}{ll}
		U & U_{\perp}
	\end{array}\right]\left[\begin{array}{cc}
		\mathbb{R}^{r \times r} & \mathbb{R}^{r \times(n-r)} \\
		\mathbb{R}^{(m-r) \times r} & 0^{(m-r) \times(n-r)} 
	\end{array}\right]\left[\begin{array}{ll}
		V & V_{\perp}
	\end{array}\right]^{T}\right\},
\end{aligned}
$$
and the normal space is 
$$
N_{X} \mathcal{M}=\left\{\left[\begin{array}{ll}
	U & U_{\perp}
\end{array}\right]\left[\begin{array}{cc}
	0^{r \times r} & 0^{r \times(n-r)} \\
	0^{(m-r) \times r} & \mathbb{R}^{(m-r) \times(n-r)}
\end{array}\right]\left[\begin{array}{ll}
	V & V_{\perp}
\end{array}\right]^{\top}\right\}.
$$
Let $P_{U}=UU^{\top}$, $P_{U}^{\perp}=I-UU^{\top}$ and  $P_{V}=VV^{\top}$, $P_{V}^{\perp}=I-VV^{\top}$. The projection to tangent space can be written as
$$
\Pi_{X}(Y)=P_{U} Y P_{V}+P_{U}^{\perp} Y P_{V}+P_{U} Y P_{V}^{\perp}.
$$
For simplicity, we rewrite (\ref{eq:rmc-full}) into the following problem:
\begin{equation}\label{eq:rmc-trans}
	\begin{array}{ll}
		\min _{X,p\in \mathbb{R}^{m \times n}} & \left\|p\right\|_{1} \\
		\text { s.t. } & P_{\Omega}(X-A) - p =0,\\
		& X \in Fr(m, n, r) .
	\end{array}
\end{equation}
The Lagrangian of (\ref{eq:rmc-trans}) can be written as $L(X,p,y) = \|p\|_1+\left\langle P_{\Omega}(X-A)-p, y\right\rangle$. It is easy to see that the KKT condition is
\begin{equation}\label{eq:rmc-kkt}
	\left\{\begin{array}{l}
		\Pi_{X}( P_{\Omega} (y) )=0, \\
		0 \in \partial\|p\|_{1}-y, \\
		P_{\Omega}(X-A)-p=0.
	\end{array}\right.
\end{equation}

By (\ref{eq:man-srcq-perturb}), the M-SRCQ is said to hold at a stationary point $X$ with respect to a multiplier $z$ if  $Dg_1(X)T_{X}\mathcal{ M}+\mathcal{C}_{\theta,g_1}(X,y)=\mathbb{R}^{m\times n} $. Since $g_1(X)=P_{\Omega}(X-A)$, we have $Dg_1(X)T_{X}\mathcal{ M}=g_1'(X)T_{X}\mathcal{ M}=P_{\Omega}(T_{X}\mathcal{ M})$. We can further obtain that $\mathcal{C}_{\theta, g_1}(X, y)=\left\lbrace d\in\mathbb{R}^{m\times n}\mid \theta^{\downarrow}(P_{\Omega}(X-A) ; d) = \left\langle d,y \right\rangle \right\rbrace $, in which
$$
\theta^{\downarrow}(P_{\Omega}(X-A) ; d)
= \sum_{P_{\Omega}(X-A)_{ij}=0}|d_{ij}|+\sum_{P_{\Omega}(X-A)_{ij}>0}d_{ij}-\sum_{P_{\Omega}(X-A)_{ij}<0}d_{ij}.
$$
Therefore, the M-SRCQ condition for RMC problem at  $(X,y)$ satisfying (\ref{eq:rmc-kkt}) is given by 
\begin{equation}\label{eq:rmc-srcq}
	P_{\Omega}(T_{X}\mathcal{ M})+\mathcal{C}_{\theta, g_1}(X, y) = \mathbb{R}^{m\times n}.
\end{equation}

The M-SOSC condition (\ref{eq:man-sosc-perturb}) can also be applied to problem (\ref{eq:rmc-trans}) since $\mathcal{ Q}=\{0\}^{m\times n}$ and the epigraph of $\theta$ are both polyhedral. It is easy to obtain that the critical cone at a stationary point $X$ is
$$
\mathcal{C}(X)=\left\lbrace \xi\in T_{X}\mathcal{M}\mid \sum_{P_{\Omega}(X-A)_{ij}=0}|\xi_{ij}|+\sum_{P_{\Omega}(X-A)_{ij}>0}\xi_{ij}-\sum_{P_{\Omega}(X-A)_{ij}<0}\xi_{ij}=0 \right\rbrace .
$$
Since 
\begin{equation}\label{eq:rmc-second-theta}
	\begin{aligned}
		&\theta_{-}^{\downarrow \downarrow}\left(g(X) ; D g(X) \xi, w \right)  = \theta_{-}^{\downarrow \downarrow}\left(P_{\Omega}(X-A) ; P_{\Omega}(\xi), w \right) \\
		= &\sum_{P_{\Omega}(X-A)_{ij}=0,\atop P_{\Omega}(\xi)_{ij}=0}|w_{ij}|+\sum_{P_{\Omega}(X-A) _{ij}>0 \text{ or}\atop P_{\Omega}(X-A) _{ij}=0,P_{\Omega}(\xi)_{ij}>0}w_{ij}-\sum_{P_{\Omega}(X-A) _{ij}<0\text{ or}\atop P_{\Omega}(X-A) _{ij}=0,P_{\Omega}(\xi)_{ij}<0}w_{ij},
	\end{aligned}
\end{equation}
the M-SOSC condition is said to hold at $X$ if for any $ \xi\in \mathcal{C}(X)\backslash\{0\}$,
\begin{equation}\label{eq:rmc-sosc}
	\begin{aligned}
		&\sup _{y\in M(X, 0,0,0)}\left\lbrace  -\psi_{\left(g_{1}(X), D g_{1}(X) \xi\right)}^{*}(y)\right\rbrace\\
		=&\sup _{y\in M(X, 0,0,0)}\left\lbrace  -\sup_{w}\left\lbrace \langle y,w\rangle -\theta_{-}^{\downarrow \downarrow}\left(g(X) ; D g(X) \xi, w \right)\right\rbrace \right\rbrace>0.
	\end{aligned}
\end{equation}

\begin{figure}[h]
	\centering
	\includegraphics[scale=0.5]{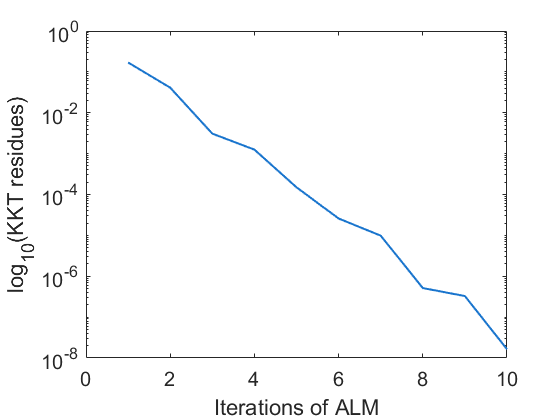}
	\caption{the KKT residues of RMC problem generated by inexact ALM }
	\label{fig:rmc-simple}
\end{figure}

We first consider a basic example of the problem (\ref{eq:rmc-full}), where $\Omega$ is the full index set.  Let $U= \begin{bmatrix}
	1& 0 &0 &0&0\\
	0 &-\dfrac{\sqrt{2}}{2}&\dfrac{\sqrt{2}}{2}&0&0 \\
	0  &\dfrac{\sqrt{2}}{2}&\dfrac{\sqrt{2}}{2}&0&0
\end{bmatrix}^{T}$, $V= \begin{bmatrix}
	1& 0 &0 &0&0\\
	0 &0.6&-0.8&0&0 \\
	0  &0.8&0.6&0&0
\end{bmatrix}^{T}$ and $S=\begin{bmatrix}
	1&0&0\\0&2&0\\0&0&3
\end{bmatrix}$. The observed matrix is set to $A=A_{\text{ex}}+E_{\text{out}}$, where $A_{\text{ex}}=USV^{T}$ is the assumed ground truth and $E_{\text{out}}$ is a matrix with random entries added only in the lower right $2\times 2$ submatrix. Since $A_{\text{ex}}$ is of rank $r=3$, $X^*=A_{\text{ex}}$ is a solution of this problem. The M-SOSC condition is satisfied since the critical cone contains only zero. Since $Dg_1(X^*)T_{X^*}F(m,n,r)=T_{X^*}F(m,n,r)$, the M-SRCQ condition holds if there exists a multiplier $y^*$ at $X^*$, such that $N_{X^*}F(m,n,r)\subset \mathcal{ C}_{\theta,g_1}(X^*,y^*).$ 
For any multiplier $y$  satisfying $Dg_1(X^*)^*y=0$, we know that $y\in N_{X^*}F(m,n,r)$, and there exists $Y\in \mathbb{R}^{m\times n}$ such that $y=P_{U}^{\perp}YP_{V}^{\perp}$.  Then for any $d=P_{U}^{\perp}\widetilde{Y}P_{V}^{\perp}\in N_{X^*}F(m,n,r)$,  if  $d\in\mathcal{ C}_{\theta,g_1}(X^*,y)$, then \begin{eqnarray*}
	\sum_{i,j\in\{4,5\}}\operatorname{sgn}(E_{\text{out}}^{ij})\widetilde{Y}_{ij}&=&\sum_{E_{\text{out}}}\operatorname{sgn}(E_{\text{out}})d_{ij}=\theta^{\downarrow}(g_1(X^*);d)=\left\langle d,y\right\rangle =\big\langle P_{U}^{\perp}\widetilde{Y}P_{V}^{\perp},P_{U}^{\perp}YP_{V}^{\perp}\big\rangle\\ &=&\sum_{i,j\in\{4,5\}}\widetilde{Y}_{ij}Y_{ij}.
\end{eqnarray*} 
It is easy to see that there exists a $y^*$ that fulfills the above equalities. Thus, the M-SRCQ condition holds at $X^*=A_{\text{ex}}$  for $y^*$. In Figure \ref{fig:rmc-simple}, we illustrate the linear convergence rate of the Riemannian ALM (Algorithm \ref{alg:alm}) in which a semismooth Newton method is employed for solving the subproblem \eqref{eq:subproblem} (see \cite[Algorithm 4.1]{ZBDZ21} for more details).

Now, look back at the original model, where $\Omega$ is no longer the full index set. In our next experiment, $A$ is generated as $A = A_{ex} + E_{out}$, where $A_{ex}$ is the ground truth and generated by $A_{ex}=L R^{\top} \in \mathbb{R}^{m \times n}$, where $L \in \mathbb{R}^{m \times r}, R \in \mathbb{R}^{n \times r}$ are two random matrices with i.i.d. standard normal distribution. $E_{out}$ is a sparse matrix with $3\%$ of the sample number are nonzero entries satisfying the exponential distribution with mean 10. The sample set is randomly generated on $\{1,2, \cdots, m\} \times \{1, 2,\cdots, n\}$ with the uniform distribution of $|\Omega| /(m n)$, and the sample size is taken as $\operatorname{OS}(m+n-r) r$, where OS is the oversampling rate for $A$ introduced in \cite{V13}, and $(m+n-r) r$ is the dimension of $Fr(m,n,r)$.

The stopping criteria is based on the KKT conditions. 
By (\ref{eq:rmc-kkt}), we terminate our ALM algorithm when the following conditions of the KKT residuals are satisfied 
$$
\max \left\lbrace \|	\Pi_{X} (P_{\Omega} (z))\|_{F},\|P_{\Omega}(X-A)-Y\|_{F},
\|Y-\operatorname{prox}_{\|\cdot\|_1}(z+Y)\|_{F}\right\rbrace \leq 1\times 10^{-7}.
$$

\begin{figure}[htbp]

	\centering
	\begin{subfigure}{0.325\linewidth}
		\centering
		\includegraphics[width=0.9\linewidth]{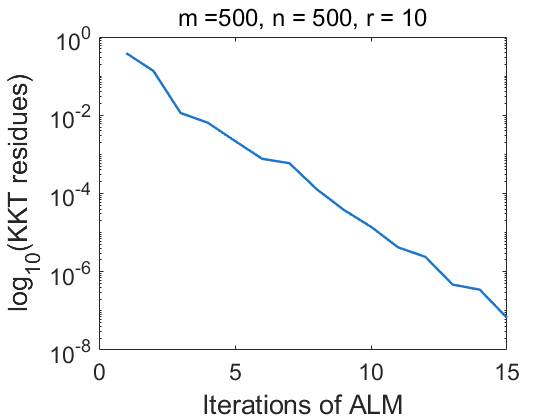}
		\caption{m=n=500, r=10}
		\label{fig:rmc500}
	\end{subfigure}
	\centering
	\begin{subfigure}{0.325\linewidth}
		\centering
		\includegraphics[width=0.9\linewidth]{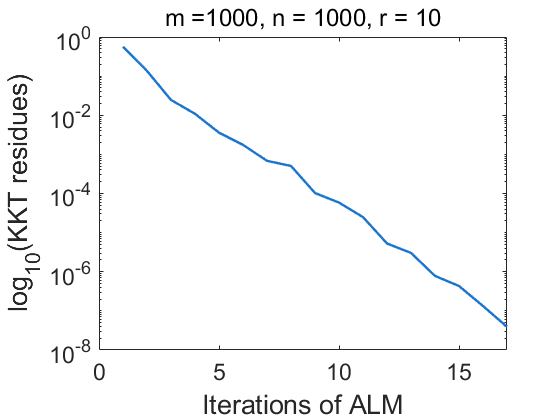}
		\caption{m=n=1000, r=10}
		\label{fig:rmc1000}
	\end{subfigure}
	\centering
	\begin{subfigure}{0.325\linewidth}
		\centering
		\includegraphics[width=0.9\linewidth]{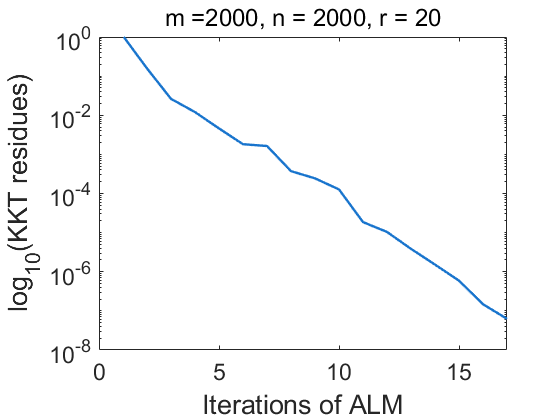}
		\caption{m=n=2000, r=20}
		\label{fig:rmc2000}
	\end{subfigure}
	\caption{The KKT residues of robust matrix completion problems generated by inexact ALM}
	\label{fig:rmcfigure}
\end{figure}

In our algorithm, a first-order method is employed to obtain a good initial point, and the semismooth Newton method is used when $\|\operatorname{grad} L_{\rho^k}\| \leq 5 \times 10^{-3}$. It should be noted that verifying the M-SRCQ and M-SOSC conditions for random examples is challenging. Despite this, Figure \ref{fig:rmcfigure} demonstrates that Riemannian ALM exhibits linear convergence when applied to these randomly generated problems. These results suggest that there might be weaker conditions under which local linear convergence of the Riemannian ALM can still be achieved.
The detailed numerical results are displayed in Table \ref{table:rmc}. The last column in Table \ref{table:rmc} shows an exact recovery of $A_{ex}$, although we do not know whether it is the global solution of problem (\ref{eq:rmc-full}).

\begin{table}[h]
	\centering
	\caption{The performance of the Riemannian augmented Lagrangian method for the robust matrix completion problem.}
	\label{table:rmc}
	\begin{tabular}{c c c|c c c c}
		\toprule
		m & n& r & iteration & time(sec) & maximum KKT residual & $\left\| X-A_{\text{ex}}\right\| $\\
		\midrule
		500 & 500 & 10 & 15 & 5.24 & 4.8339e-08 & 1.6186e-08\\
		\midrule
		1000 & 1000 & 10 & 17 & 14.87 &3.5088e-08 & 8.3608e-08\\
		\midrule
		2000 & 2000 & 20 & 17 & 33.47 &2.5007e-08 & 2.4877e-08\\
		\midrule
		5000 & 5000 & 20 & 30 & 430.02 &6.6165e-09 & 3.4093e-09\\
		\bottomrule
	\end{tabular}
\end{table}

\section{Conclusion}\label{sec:conslusion}

This paper studies the strict Robinson constraint qualification and the second order condition for the nonsmooth optimization problems on manifolds. We show that the  M-SRCQ and M-SOSC conditions are equivalent to the robust isolated calmness of the  KKT solution mapping. Under these two conditions, we show that the iterations generated by the Riemannian augmented Lagrangian method converge to a KKT point, and the rate of convergence is linear. Numerical results on a class of nonsmooth optimizations over the sphere and the robust matrix completion problems demonstrate the convergence rate, respectively. 

The work done in this paper on the perturbation properties of the nonsmooth manifold optimization problem (\ref{eq:prime}) and the local convergence analysis of the Riemannian augmented Lagrangian method is by no means complete. Due to the rapid advances in the applications of nonsmooth optimization problems on manifolds in different fields, the study of manifold optimization problems will become even more important, and many other properties of perturbation analysis and algorithm design are waiting to be explored. For example, consider the sparse principal component analysis (SPCA) problem (cf. e.g., \cite{JTU03,ZHT06})
\begin{equation}\label{eq:spca}
	\begin{array}{ll}
		\min _{X \in \mathbb{R}^{n \times p}} & -\operatorname{tr}\left(X^{\top} A^{\top} A X\right)+\mu\|X\|_{1} \\
		\text { s.t. } & X^{\top} X=I_{p},
	\end{array}
\end{equation}
where $\operatorname{St}(n,p)=\{X\in \mathbb{R}^{n\times p}\mid X^{\top}X=I_p\}$ and $f:\operatorname{St}(n,p)\rightarrow \mathbb{R}$ is a smooth function.
It is worth noting that unlike the nonsmooth optimization problem on sphere \eqref{eq:sphere},  the local optimal solutions of the SPCA problem \eqref{eq:spca} may not satisfy the M-SRCQ condition if $p>1$. On the other hand, the numerical experiments given in \cite{ZBDZ21} applying the Riemannian ALM to SPCA problems have shown that the Algorithm \ref{alg:alm} performs well numerically and usually can be observed a local fast linear convergence rate (cf. \cite{ZBDZ21}). Therefore, the weaker condition for ensuring the local linear convergence rate of the Riemannian ALM of nonsmooth optimization problems on manifolds is an interesting and important issue for our future work.  Another critical and practical issue is designing efficient algorithms for solving the Riemannian augmented Lagrangian subproblem \eqref{eq:subproblem} in Algorithm \ref{alg:alm}.

\bibliographystyle{spmpsci}
\bibliography{almtheta}

\end{document}